\RequirePackage{fix-cm}
\documentclass{svjour3}                     
\smartqed  

\usepackage{amsmath,amsfonts,amssymb,mathrsfs,cases,mathdots,color,colordvi,url,graphicx,multirow,caption,subcaption,float,lscape,algorithm,algorithmic}
\usepackage{geometry}
\usepackage{hyperref}
\usepackage{nicefrac}
\usepackage[justification=centering]{caption}
\usepackage[title]{appendix}
\newtheorem{assumption}{Assumption}

\newcommand{\BOX}{\hfill \Box}

\begin{document}

\title{Local convergence analysis of augmented Lagrangian method for nonlinear semidefinite programming}

\author{Shiwei Wang \and
Chao Ding
}
\institute{S.W. Wang \at
      Institute of Applied Mathematics, Academy of Mathematics and Systems Science, Chinese Academy of Sciences, Beijing,  P.R. China,  School of Mathematical Sciences, University of Chinese Academy of Science, Beijing,  P.R. China. \\
       \email{wangshiwei182@mails.ucas.ac.cn}%
	\and
	C. Ding \at
	Institute of Applied Mathematics, Academy of Mathematics and Systems Science, Chinese Academy of Sciences, Beijing,  P.R. China. The research of this author was supported by the National Natural Science Foundation of China under project No. 12071464 and the Beijing Natural Science Foundation (Z190002). \\
	\email{dingchao@amss.ac.cn}
}

\date{This version: Aug 8, 2022}

\maketitle
\begin{abstract}
The augmented Lagrangian method (ALM) has gained tremendous popularity for its elegant theory and impressive numerical performance since it was proposed by Hestenes and Powell in 1969. It has been widely used in numerous efficient solvers to improve numerical performance to solve many problems. In this paper, without requiring the uniqueness of multipliers, the local (asymptotic Q-superlinear) Q-linear convergence rate of the primal-dual sequences generated by ALM for the nonlinear semidefinite programming (NLSDP) is established by assuming the second-order sufficient condition (SOSC) and the semi-isolated calmness of the Karush–Kuhn–Tucker (KKT) solution under some mild conditions.


\end{abstract}

\section{Introduction}\label{sec1}

The non-convex semidefinite programming problem is attracting more attention for its wide applications in machine learning, structural design, and other fields. In the well-known library COMPleib \cite{Leibfritz05}, there are about 168 test examples for nonlinear semidefinite programs, control system design, and related problems. In this paper, we consider the nonlinear semidefinite programming (NLSDP) problem in the following form:
\begin{equation}\label{eq:NLSDP}
	\begin{array}{cl}
		\displaystyle{\min_{x\in {\cal X}}} & f(x) \\ [3pt]
		{\rm s.t.} & h(x)=0, \\ [3pt]
		&{G(x)\in {\cal S}^n_+},
	\end{array} 
\end{equation}
where 
$f:{\cal X}\to \Re$, $h:{\cal X}\to {\cal Y}$ and $G:{\cal X}\to {\cal S}^n$ are twice continuously differentiable, ${\cal X}$, ${\cal Y}$ are given Euclidean spaces, ${\cal S}^n$ is the linear space of all $n\times n$ real symmetric matrices equipped with the usual Frobenius inner product and its induced norm. For notational simplicity, we use $\langle \cdot, \cdot\rangle $ to denote inner product of every Euclidean spaces and $\|\cdot\|$ to denote its induced norm. The exact meaning of these notations can be deduced from the context. ${\cal S}^n_+$ (${\cal S}^n_-$) is used to represent the $n$-dimensional positive (negative) semidefinite cone. 
For a general constraint optimization problem, if the constraint set ${\cal K}$ is some polyhedron, we say the problem is polyhedral; otherwise, we say the problem is non-polyhedral. Especially, when ${\cal K}=\{0\}\times\Re^n_+$, it is the well-known nonlinear programming (NLP). It is easy to see that NLSDP is a non-polyhedral problem as it possesses the positive semidefinite (SDP) cone constraint. 

In this paper, we will mainly focus on the local convergence analysis of augmented Lagrangian method (ALM)  for \eqref{eq:NLSDP}. The Lagrangian function of problem \eqref{eq:NLSDP}  
is defined by 
\begin{equation}\label{eq:L-function}
	L(x,y,\Gamma):=f(x)+\langle y,h(x)\rangle  +\langle \Gamma,G(x)\rangle, \quad (x,y,\Gamma)\in{\cal X}\times\mathcal{Y}\times{\cal S}^n.
\end{equation}
For any $(y,\Gamma)\in\mathcal{Y}\times{\cal S}^n$, denote the first-order and second-order derivatives of $L(\cdot,y,\Gamma)$ at $x\in{\cal X}$ by $\nabla_xL(x,y,\Gamma)$ and $\nabla_{xx}^2L(x,y,\Gamma)$, respectively. 
Augmented Lagrangian function is firstly introduced by Arrow and Solow in 1958 \cite{ASolow58} to study a differential equation method for solving equality constrained optimization problems. It is identified by Rockafellar in 1970 \cite{Rockafellar70} and firstly studied in detail by Buys in his doctoral dissertation in 1973 \cite{Buys72} to inequality constraints. 
The augmented Lagrangian function of \eqref{eq:NLSDP} takes the following form (cf. \cite[Section 11.K]{RoWe98} and \cite{ShapS04})
\begin{equation}\label{eq:deflag}
	\mathscr{L}(x,\zeta,\rho):=f(x)+\frac{\rho}{2}{\rm dist}^2(\Phi(x)+\frac{\zeta}{\rho},{\cal K})-\frac{\|\zeta\|^2}{2\rho},
\end{equation}
where $\zeta:=(y,\Gamma)$, $\Phi(x):=(h(x),G(x))$, ${\cal K}:=\{0\}\times{\cal S}_+^n$. And for a given $z$, ${\rm dist}(z,{\cal K})$ is the distance from point $z$ to ${\cal K}$.  
For a given initial point $(x^0,\zeta^0)\in{\cal X}\times\mathcal{Y}\times{\cal S}^n$ and a constant $\rho^0>0$, the $(k+1)$-th iteration of ALM for NLSDP \eqref{eq:NLSDP} is 
\begin{equation*}
	\left\{\begin{array}{ll}
		x^{k+1}\approx\arg\min\{\mathscr{L}(x,\zeta^{k},\rho^k)\},  \label{eq:subp1}\\
		\zeta^{k+1}=\rho^k\big[\Phi(x^{k+1})+\frac{\zeta^k}{\rho^k}-\Pi_{\cal K}(\Phi(x^{k+1})+\frac{\zeta^k}{\rho^k})\big]  
	\end{array}\right.
\end{equation*}
with $\rho^{k+1}$ updated by certain rule.

The augmented Lagrangian method was firstly proposed by Hestenes \cite{Hes69} and Powell \cite{Po72} for solving the equality constrained problem and was generalized by Rockafellar \cite{Rockafellar1973} to NLP. It rapidly grew into popularity for its mathematical elegance and impressive numerical performance in various areas, like statistical optimization (e.g., Lasso \cite{Tibshirani96}), machine learning and game theory. It has also been implemented in many powerful large scale solvers like SDPNAL$+$ \cite{YSToh15,ZSToh10}, QSDPNAL \cite{LSToh18}, SuitedLasso \cite{LSToh182} and so on. 

There are also many works focused on the theory of this algorithm. For convergence analysis of ALM, tremendous work has been established since it was proposed. 
Powell \cite{Po72} demonstrated that for the equality constrained problem, if the second-order sufficient condition (SOSC) and linear independence constraint qualification (LICQ) were satisfied, the algorithm should converge locally at a linear rate, without the need for having $\rho\rightarrow\infty$. This implies that ALM may provide numerical stability, which the usual penalty methods do not possess. 
In 1973, Rockafellar \cite{Rockafellarhp} and Tretykov \cite{Tretyakov73} proved the global convergence of the augmented Lagrangian method for convex optimization problem with inequality constraints for any $\rho>0$ based on the saddle point theorem established in \cite{Rockafellar1973}.

For the convex NLP problem, the local convergence rate of the ALM can be derived through its deep connection with the dual proximal point algorithm (PPA) as studied by Rockafellar in \cite{Rockafellar1976}. As stated in \cite[Proposition 3, Theorem 2]{Rockafellar1976}, one can obtain the Q-linear convergence rate of the dual sequence generated by the ALM under the upper Lipschitz continuity of the dual solution mapping at the origin, the boundedness of dual sequence and certain stopping criteria on the inexact computations of the augmented Lagrangian subproblems. For more details about PPA and monotone operators, please see \cite{Rockafellarmo,Rockafellar1976,Pennane02}.

Following this way, the convergence rate of ALM for general convex optimization problems can also be attained under very mild conditions with implementable stopping criteria for the ALM subproblems.
In 1984, Luque relaxed the upper Lipschitz continuity of the dual solution mapping used in \cite{Rockafellar1976}, which required the uniqueness of the optimal solution, by an error bound type condition \cite[(2.1)]{Luque} that is known to be satisfied for polyhedron \cite{Robinson81} but difficult to be verified for non-polyhedron. 
 In 2019, Cui et al. \cite{CSToh2016} established the asymptotic R-superlinear convergence of the KKT residuals and asymptotic Q-superlinear convergence of the dual sequence generated by the ALM for solving convex NLSDP, under a quadratic growth condition on the dual problem that neither local solution nor the multiplier is required to be unique. 
 Their remarkable work improved \cite{Rockafellar1976} in giving a practical stopping criterion for ALM subproblem under the Robinson constraint qualification (RCQ) (for the improvement of implementable stopping criteria, see also \cite{ESilva13}) and obtaining the convergence of the KKT residuals 
 with the application of KKT residual information. Also, they relaxed Luque's condition by the calmness of the dual solution mapping at the origin.

When it comes to non-convex optimization problems, fruitful results have been established for the polyhedral case. In 1982, Bertsekas \cite{Bertsekas} established that the generated dual sequence converges Q-linearly and the corresponding primal sequence converges R-linearly under SOSC, LICQ and the strict complementarity for NLP. His result shows that the ratio constant is proportional to $1/\rho$, which implies the convergence can be accelerated by increasing $\rho$. Efforts are made to weaken the above conditions. Firstly, successful attempts are made to remove the strict complementarity condition, e.g., Conn et al. \cite{Conn91}, Contesse-Becker \cite{Lcb93}, and Ito and Kunisch \cite{IKunisch90} derived linear convergence rate for the ALM of general NLP. Secondly, it is also crucial to weaken LICQ, which implies the uniqueness of multipliers. As in real-world, multipliers are usually non-unique, e.g., considering Lasso as a dual problem. In 2012, Fernandez and Solodov \cite{FSolodov12} firstly studied this topic for NLP without requiring the multiplier to be unique. This work is a milestone to establishing the convergence by removing the uniqueness of the Lagrangian multiplier and strict complementarity. Recently Hang and Sarabi \cite{HSarabi20} established the local convergence for piecewise linear quadratic composite optimization problems under merely SOSC, which inspires us to study whether their results can be extended to NLSDP \eqref{eq:NLSDP}. Their success relies on the validity of upper Lipschitz continuous of KKT solution mapping when SOSC is satisfied, see \cite{DRoc97,IKur13,Klatte00,MSarabi2018,HSarabi20}. However, this does not hold for non-polyhedral case as mentioned in \cite{CSToh2016} by using \cite[Example 4.54]{bonnans}. For comprehensive surveys about the augmented Lagrangian method for nonlinear programming, see \cite{Bertsekas,GTre89,Rocksurvey}. 
Recently, Rockafellar \cite{Rockafellarcal} shaded lights on how to derive ALM convergence rate for non-convex ``fully amenable" \cite[10F]{RoWe98} problems through its connection with PPA for the dual problem and the variational sufficiency studied in \cite[Page 6]{Rockafellarc}. He successfully obtained the ALM primal R-linear convergence from the ALM dual Q-linear convergence for generalized NLP by assuming variational sufficiency. However, NLSDP does not belong to the fully amenable kind. Also, as mentioned in \cite{Rockafellarcal}, variational sufficiency may fail even under SOSC. 

For non-convex non-polyhedral problem, Sun et al. \cite{SSZhang08} proved the convergence rate of NLSDP under strongly SOSC \cite{Sun06} together with nondegeneracy (cf. \cite{SSZhang08} or \cite{bonnans}). In 2019, Kanzow and Steck \cite{KSteck19} justified the primal-dual linear convergence of ALM under SOSC and strong Robinson constraint qualification (SRCQ) for $C^2$-cone reducible constrained problems, which include NLSDP and nonlinear second-order cone programming (NLSOC). Recently, the primal-dual linear convergence rate of ALM for NLSOC  is also studied under SOSC and the semi-isolated calmness of the KKT solution mapping (see Definition \ref{def:semi-ic}) in \cite{HMSarabi} with the multiple uniqueness assumption. Going through all the papers above, it is not hard to see that existing works usually suppose either the problem is convex (or polyhedral), or the Lagrangian multiplier is (locally) unique. However, few results on the local convergence rate of ALM have been established for non-convex non-polyhedral problems without the multiple uniqueness.

In this paper, under some mild conditions, we establish a locally (asymptotic Q-superlinear) Q-linear convergence of the primal-dual sequences generated by ALM without assuming the uniqueness of Lagrangian multipliers for NLSDP under SOSC and the semi-isolated calmness of the KKT solution mapping. Furthermore, for NLSDP, we provide a sufficient condition for the semi-isolated calmness of the KKT solution mapping. The remaining parts of this paper are organized as follows. In the next section, we introduce some preliminary knowledge in semidefinite cone and variational analysis. In Section \ref{sec:usoe}, we study the uniform second-order expansion for the Moreau envelop of the indicator function of the SDP cone, which is essential to the main convergence result. In Section \ref{sec5}, we obtain the main result on the local (asymptotic Q-superlinear) Q-linear convergence of ALM for NLSDP with the help of uniform second-order growth condition of augmented Lagrangian function. Section \ref{sec:scsic} is devoted to the sufficient condition of semi-isolated calmness of the set of KKT points. Also in this section, we illustrate by two examples that the conditions proposed in our main result can be satisfied. We conclude our paper and make some comments in the final section.

\section{Preliminaries}
This section lists some preliminaries on the positive semidefinite cone and variational analysis,  which will be used in this paper. Detailed discussions on these subjects can be found in \cite{Clarke,Mor,RoWe98}.

Let  $A \in {\cal S}^n$ be given. We use $\lambda_{1}(A)\ge \lambda_2(A) \ge \ldots \ge \lambda_{n}(A)$ to denote the eigenvalues of $A$ (all real and counting multiplicity) arranging in  nonincreasing order and use $\lambda(A)$ to denote the vector of the ordered eigenvalues of $A$. Let $\Lambda(A):= {\rm Diag}(\lambda(A))$. Also, we use $v_1(A)> \dots> v_d(A)$ to denote the different eigenvalues. Consider the eigenvalue decomposition of $A$, i.e., $A={P} \Lambda(A){P}^{T}$, where ${P}\in{\cal O}^{n}(A)$ is a corresponding
orthogonal matrix of the orthonormal eigenvectors. 
By considering the index sets of positive, zero, and negative eigenvalues of $A$, we are able to write $A$ in the following form
\begin{equation}\label{eq:eig-decomp}
	A= \left[\begin{array}{ccc}
		{P}_{\alpha} & {P}_{\beta} & {P}_{\gamma}
	\end{array}\right] \left[\begin{array}{ccc}
		\Lambda(A)_{\alpha\alpha} & 0 & 0 \\ [3pt]
		0 & 0 & 0\\ [3pt]
		0 & 0 & \Lambda(A)_{\gamma\gamma}
	\end{array}\right]\left[\begin{array}{c}
		{P}_{\alpha}^T \\ [3pt]
		{P}_{\beta}^T \\
		{P}_{\gamma}^T
	\end{array}\right].
\end{equation}
where $\alpha:=\{i: \lambda_i(A)>0\}$, $\beta:=\{i: \lambda_i(A)=0\}$ and $\gamma:=\{i: \lambda_i(A)<0\}$. 
It is known from \cite[Example 3.140]{bonnans} that the $n$-dimensional positive semidefinite cone ${\cal S}_+^n$ is $C^{\infty}$-cone reducible. 
We use $\Pi_{{\cal S}_+^n}(A)$ to represent the projection from $A$ to ${\cal S}_+^n$ and it follows from \cite[p.5]{Sun06} that 
\begin{equation*}
	\Pi_{{\cal S}_+^n}(A)= \left[\begin{array}{ccc}
		{P}_{\alpha} & {P}_{\beta} & {P}_{\gamma}
	\end{array}\right] \left[\begin{array}{ccc}
		\Lambda(A)_{\alpha\alpha} & 0 & 0 \\ [3pt]
		0 & 0 & 0\\ [3pt]
		0 & 0 & 0
	\end{array}\right]\left[\begin{array}{c}
		{P}_{\alpha}^T \\ [3pt]
		{P}_{\beta}^T \\
		{P}_{\gamma}^T
	\end{array}\right].
\end{equation*}
From \cite[Theorem 4.7]{SSun02} we know that the metric projection operator $\Pi_{{\cal S}_+^n}(\cdot)$ is directionally differentiable at any $A\in{\cal S}^n$ and the directional derivative of $\Pi_{{\cal S}_+^n}(\cdot)$ at $A$ along direction $H\in{\cal S}^n$ is given by 
\begin{equation} \label{eq:dd-projection}
	\Pi_{{\cal S}_+^n} ^{\prime}(A;H) =  {P} \left[
	\begin{array}{ccc}
		{\widetilde H}_ {\alpha\alpha}   & {\widetilde H}_{\alpha\beta} &\Sigma_{\alpha\gamma} \circ {\widetilde H}_{\alpha\gamma}\\
		\widetilde{H}_{\alpha\beta}^T &  \Pi_{{\cal S}_+^{\rvert\beta\lvert} }({\widetilde H}_{\beta\beta}) &0\\
		{\Sigma}_{\alpha\gamma}^T \circ {\widetilde H}_{\alpha\gamma}^T & 0 & 0
	\end{array}
	\right]  {P}^T,
\end{equation}
where ${\widetilde H}:={P}^T H{P}$, ``$\circ$" is the Hadamard product and
\begin{equation}\label{eq:def-Sigma}
	\Sigma_{ij}:=\frac{\max\{\lambda_{i}(A),0\}-\max\{\lambda_{j}(A),0\}}{\lambda_{i}(A)-\lambda_{j}(A)},\quad i,j=1,\ldots,n,
\end{equation} where $0/0$ is defined to be $1$. Moreover, we also have 
\begin{equation}\label{eq:dd-projection-}
\Pi_{{\cal S}_-^n} ^{\prime}(A;H)=H-\Pi_{{\cal S}_+^n} ^{\prime}(A;H)=P\left[
	\begin{array}{ccc}
		0   & 0 &(E-\Sigma)_{\alpha\gamma} \circ {\widetilde H}_{\alpha\gamma}\\
		0 &  \Pi_{{\cal S}_-^{\rvert\beta\lvert} }({\widetilde H}_{\beta\beta}) &{\widetilde H}_{\beta\gamma}\\
		(E-{\Sigma})_{\alpha\gamma}^T \circ {\widetilde H}_{\alpha\gamma}^T & {\widetilde H}_{\beta\gamma}^T & {\widetilde H}_{\gamma\gamma}
	\end{array}
	\right]P^T.
\end{equation} 
For more details on the properties of SDP cone, we refer the reader to \cite{Sun06,SSZhang08,DSYe14}.

For a given Euclidean space ${\cal X}$. 
Let $C$ be any subset in ${\cal X}$, the Bouligand tangent/contingent cone of $C$ at $x$ is a closed cone defined by
\[
T_{C}(x):= \left\{d\in{\cal X}\mid \mbox{$\exists\, t^k\downarrow 0$ and  $d^k\to d$ with  $x+t^kd^k\in C$ for all $k$} \right\}.
\]
The regular/Fr\'{e}chet normal cone of $C$ at $x$ is defined by
\[
\widehat{N}_C(x):=\left\{v\in {\cal X}\mid \langle v,x'-x\rangle \le o(\|x'-x\|)\ \forall\,x'\in C  \right\}.
\]
The limiting/Mordukhovich normal cone is defined by
\[
N_C(x):=\left\{\lim_{k\to\infty}v^k\mid v^k\in \widehat{N}_{C}(x^k),\ x^k\to x,\  \mbox{$x^k\in C$ for all $k$}  \right\}.
\]
When $C$ is convex, the regular normal cone and limiting normal cone coincide with the normal cone in the sense of  convex analysis \cite{Rockafellar1970}, i.e.,
${N}_{C}(x)=\widehat{N}_C(x)=\left \{v \in \mathcal{X}:\langle v,x'-x\rangle\leq 0 \ \forall\,x'\in C \right\}$.  
For any $x\in C$, the critical cone  associated with ${ y}\in{N}_C({x})$ of $C$ is defined in \cite[Page 98]{DoRo14}, i.e., ${\cal C}_{C}({x},{y})=T_{C}({x})\cap({y})^{\perp}$. It is well known (see e.g., \cite[(19)]{Sun06}) that the critical cone of SDP at a given $Y\in{ N}_{{\cal S}_+^n}(X)$ can be completely described as
\begin{equation}\label{eq:crisdp}
{\cal C}_{{\cal S}_+^n}(X,Y):=\{U\in{\cal S}^n\mid {P}^T_{\beta}U{P}_{\beta}\in{\cal S}_+^{\lvert\beta\rvert},\; {P}^T_{\beta}U{P}_{\gamma}=0,\;{P}^T_{\gamma}U{P}_{\gamma}=0\},
\end{equation}
where $X+Y$ has the eigenvalue decomposition in \eqref{eq:eig-decomp}.

Given a proper, lower semi-continuous function $f:{\cal X}\rightarrow(-\infty,+\infty]$ and a constant $\rho>0$, the Moreau envelop of function $f$ at $y\in{\cal X}$  
is defined by
\begin{equation}\label{eqdef:more}
	e_{\rho}f({y}):=\inf\limits_{{x}\in\mathcal{X}}\left\{f({ x})+\frac{1}{2\rho}\|{ x}-{ y}\|^2\right\}. 
\end{equation}
In particular, when $\rho=1$, we denote $ef({y)}$ as $ef({y})$ for simplicity. 
Especially, when $f$ is the indicator function of ${\cal S}_+^n$ with $\rho=1$ at $A\in{\cal S}^n$, we denote it as $e\delta_{{\cal S}_+^n}(A)$. 

The following definition of second-order subderivative is taken from \cite[Definition 13.3]{RoWe98}. 
\begin{definition}\label{eqdef:sosd}
	For $f:{\cal X}\rightarrow(-\infty,+\infty]$, any $\bar{x}\in{\cal X}$ with $f(\bar{x})$ finite and any $\bar{y}\in{\cal X}$, the second subderivative of $f$ at $\bar{x}$ for $\bar{y}$ is defined as 
	\begin{equation*}
		{\rm d}^2f(\bar{x},\bar{y})(v)=\liminf_{c\downarrow0,v'\rightarrow v}\frac{f(\bar{x}+cv')-f(\bar{x})-c\langle\bar{y},v\rangle}{\frac{1}{2}c^2}. 
	\end{equation*}
	The second semiderivative of $f$ at $\bar{x}$ is 
	\begin{equation*}
		{\rm d}^2f(\bar{x})(v)=\lim_{c\downarrow0,v'\rightarrow v}\frac{f(\bar{x}+cv')-f(\bar{x})-c\,{\rm d}f(\bar{x})(v)}{\frac{1}{2}c^2}, 
	\end{equation*}
	where ${\rm d}f(\bar{x})(v)$ is the semiderivative of $f$ at $\bar{x}$ defined in \cite[Definition 7.20]{RoWe98}. 
\end{definition}

Especially, for augmented Lagrangian function of NLSDP \eqref{eq:deflag}, we use ${\rm d}^2_x\mathscr{L}$ to denote the partial second semiderivative on $x$. The following definition  is an extension of the second order expansion of functions defined in \cite[Definition 1.1]{PRock}.  
\begin{definition}
	Consider a function $f:{\cal X}\rightarrow\Re$ and a point $\bar{x}$ where $f$ is differentiable. We say $f$ has a uniform second order expansion at $\bar{x}$ with certain conditions if $f$ satisfies the following two conditions. 
	\begin{itemize}
		\item[(i)] $f$ has a second order expansion at $\bar{x}$, i.e., there exists a finite, continuous and positively homogeneous of degree 2 function $g$ such that 
		$$f(\bar{x}+cv)=f(\bar{x})+c\langle\nabla f(\bar{x}),v\rangle+\frac{c^2}{2}g(v)+o(c^2\|v\|^2),\quad c\in\Re,\;v\in{\cal X}.$$
		\item[(ii)] There exists a constant $r>0$ such that $o(c^2\|v\|^2)$ is uniform for all $x\in\mathbb{B}_r(\bar{x})$ with certain conditions, i.e., for all $\varepsilon>0$, there exist  positive constants $\omega$, $r$  such that for all $x\in\mathbb{B}_{r}(\bar{x})$ with certain conditions and all $\|cv\|\leq\omega$, we have 
		$$\frac{f({x}+cv)-f({x})-c\langle\nabla f({x}),v\rangle-\frac{c^2}{2}g(v)}{\|c^2v^2\|}\leq \varepsilon,$$
		where $\omega$ is uniform for all $x\in\mathbb{B}_{r}(\bar{x})$ with certain conditions. 
	\end{itemize}
\end{definition}

\section{The uniform second order expansion of $e\delta_{{\cal S}_+^n}(A)$}\label{sec:usoe}
In this section, we shall establish the uniform expansion of the Moreau envelop $e\delta_{{\cal S}_+^n}(\cdot)$ of the indicator function of ${\cal S}^n_+$, which is crucial for the subsequent analysis of deriving the uniform quadratic growth condition of augmented Lagrangian function. Before we step further, we need to give the following notation. Given two Euclidean spaces ${\cal X}$, ${\cal Z}$, a positive constant $r$ and $\overline{A}\in{\cal X}$. Considering a mapping $\Delta:{\cal X}\times{\cal X}\rightarrow{\cal Z}$. We say $\Delta(A,H)={O}(\|H\|)$ with ${O}(\|H\|)$ uniform for all $A\in\mathbb{B}_{r}(\overline{A})$ with certain conditions if there exist  positive constants $q$, $\omega$, $r$  such that for all $A\in\mathbb{B}_{r}(\overline{A})$ and all $\|H\|\leq\omega$, we have 
$$\frac{\|\Delta(A,H)\|}{\|H\|}\leq q,$$
where $\omega$ and $q$ are uniform for all $A\in\mathbb{B}_{r}(\overline{A})$ with certain conditions. We call $\omega$ the uniform radius and $q$ the uniform constant of $O(\|H\|)$. To obtain the main result of this section, firstly we need the following lemma, which illustrates the uniform expansion for eigenvalue vector matrix. Its non-uniform form was stated in \cite{SSun03} and essentially proved in the derivation of \cite[Lemma 4.12]{SSun02}. 
\begin{lemma}\label{lemma:unih}
Given a fixed $\overline{A}\in{\cal S}^n$. Let $0<r<\min_{i<j}\{ v_i(\overline{A})-v_j(\overline{A})\}/3$. For any $H\in{\cal S}^n$ and $A\in\mathbb{B}_{r}(\overline{A})$, let $U$ be an orthogonal matrix such that 
\begin{equation}\label{eq:unih2}
U^T\big(\Lambda({A})+H\big)U=\Lambda\big(\Lambda(A)+H\big).
\end{equation}
Then, for any $H\rightarrow0$, we have 
\begin{equation}\label{eq:unih1}
\left\{\begin{array}{ll}
U_{\bar{\iota}_{s} \bar{\iota}_{t}}={O}(\|H\|), & s, t=1, \cdots, \bar{d},\ s \neq t \\
U_{\bar{\iota}_{s} \bar{\iota}_{s}} U_{\bar{\iota}_{s} \bar{\iota}_{s}}^{T}=I_{\lvert\bar{\iota}_{s}\rvert}+{O}(\|H\|^{2}), & s=1, \cdots, \bar{d},
\end{array}\right.
\end{equation}
where $\bar{\iota}_s:=\iota_s(\overline{A})=\{i\mid \lambda_i(\overline{A})=v_s(\overline{A})\}$, $s=1,\dots,\bar{d}$. 
Furthermore, for each $s\in\{1,\dots,\bar{d}\}$, there exists $Q_s\in{\cal O}^{\lvert\bar{\iota}_s\rvert}$ such that 
\begin{equation}\label{eq:unih3}
U_{\bar{\iota}_{s} \bar{\iota}_{s}}=Q_k+{O}(\|H\|^2)
\end{equation}
and 
\begin{equation}\label{eq:unih4}
Q_{s}^{T} H_{\bar{\iota}_{s} \bar{\iota}_{s}} Q_{s}=\Lambda_{\bar{\iota}_{s} \bar{\iota}_{s}}(\Lambda(X)+H)-Q_s^T\Lambda(A)_{\bar{\iota}_{s} \bar{\iota}_{s}}Q_s+{O}(\|H\|^{2}).
\end{equation}
It is worth to note that the ${O}(\|H\|)$ and ${O}(\|H\|^{2})$ above are uniform for all $A\in\mathbb{B}_{r}(\overline{A})$. 
\end{lemma} 
{\bf Proof.} See Appendix \ref{app:A}. 
\vskip 10pt
By applying Lemma \ref{lemma:unih}, we can obtain the following result. 
\begin{lemma}\label{lemma:unisec}
	Given $\overline{A}\in{\cal S}^n$ and let $0<r<\min_{i<j}\{v_i(\overline{A})-v_j(\overline{A})\}/3$. For any $H\in{\cal S}^n$ and $A\in\mathbb{B}_{r}(\overline{A})$, let $U$ be an orthogonal matrix such that 
	\begin{equation}\label{eq:lemap1}
		U^T\big({A}+H\big)U=\Lambda({A}+H).
	\end{equation}
	For all $t\in\{1,\dots,\bar{d}\}$, there exist $Q_t\in{\cal O}^{\lvert\bar{\iota}_t\rvert}$ (depends on $H$) such that for all $H\rightarrow0$, 
	\begin{align*}
		&(P^TU)_{\bar{\iota}_s\bar{\iota}_t}=\Theta^{st}\circ(\widetilde{H}_{\bar{\iota}_s\bar{\iota}_t}Q_t)+{O}(\|H\|^2),\quad s\neq t,
	\end{align*}
	where ${O}(\|H\|^{2})$ is uniform for all $A\in\mathbb{B}_r(\overline{A})$, $(\Theta^{st})_{ij}=1/\big((\Lambda(A)_{\bar{\iota}_t\bar{\iota}_t})_{ii}-(\Lambda(A)_{\bar{\iota}_s\bar{\iota}_s})_{jj}\big)$ and $\widetilde{H}=P^THP$, $P\in{\cal O}^n(A)$.  
\end{lemma}
{\bf Proof.} See Appendix \ref{app:C}. 
\vskip 10pt
For two matrices $A,\overline{A}\in{\cal S}^n$, notation $\pi(A)=\pi(\overline{A})$ means these two matrices possess the same index sets of different eigenvalues, i.e, $A$ and $\overline{A}$ both have $\bar{d}$ different eigenvalues with $\iota_l(A)=\iota_l(\overline{A}):=\{i\mid \lambda_i(\overline{A})=v_s(\overline{A})\}$ for all $l=1,\dots,\bar{d}$ and the number of positive, zero, negative eigenvalues equal. 
Applying Lemma \ref{lemma:unih}, we can get the uniform 1-order B-differentiability of projection function for SDP case, which is an enhancement of \cite[Proposition 2.6]{DSYe14}. 
\begin{proposition}\label{prop:unisocpi}
	Given a fixed $\overline{A}\in{\cal S}^n$ and let $0<r<\min_{i<j}\{v_i(\overline{A})-v_j(\overline{A})\}/3$. The metric projection operator is uniformly 1-order B-differentiable for any $A\in\mathbb{B}_r(\overline{A})$ with $\pi(\overline{A})=\pi(A)$, i.e., for ${\cal S}^n\ni H\rightarrow0$, 
	\begin{equation}\label{eq:piuni}
		\Pi_{\mathcal{S}_{+}^{n}}(A+H)-\Pi_{\mathcal{S}_{+}^{n}}(A)-\Pi_{\mathcal{S}_{+}^{n}}^{\prime}(A ; H)={O}(\|H\|^{2})
	\end{equation}
	and ${O}(\|H\|^{2})$ is uniform for all $A\in\mathbb{B}_r(\overline{A})$ with $\pi(\overline{A})=\pi(A)$. 
\end{proposition}
{\bf Proof.} See Appendix \ref{app:B}. 
\vskip 10pt
With the help of the properties of the second subderivative, we can calculate its corresponding Moreau envelop in the following explicit form. This will be of great use in deriving the main result of this section. 
\begin{lemma}\label{lemma:morcal}
	Suppose  ${\cal K}=\{0\}\times{\cal S}^n_+$. Given $(\Phi({x}),\zeta)\in{\rm gph}\,N_{{\cal K}}$. Denote $A=G(x)+\Gamma$ and $A$ possesses the eigenvalue decomposition in \eqref{eq:eig-decomp}. We have that for all $H\in{\cal S}^n$, the Moreau envelop of ${\rm d}^2\delta_{{\cal K}}(\Phi({x}),\zeta)(\cdot)$ can be calculated in the following form, i.e., for all $(b,B)\in{\cal Y}\times{\cal S}^n$,  
	\begin{align*}
		&e_{1/(2\rho)}\big({\rm d}^2\delta_{{\cal K}}(\Phi({x}),\zeta)\big)(b,B)=e_{1/(2\rho)}\big({\rm d}^2\delta_{{\cal S}_+^n}(G({x}),\Gamma)\big)(B)+\rho\|b\|^2\nonumber\\
		=&2\rho\sum\limits_{i\in\beta,\,j\in\gamma}\widetilde{B}_{ij}^2+\rho\sum\limits_{i,\,j\in\gamma}\widetilde{B}_{ij}^2+\rho{\rm dist}^2(\widetilde{B}_{\beta\beta},{\cal S}_+^{\lvert\beta\rvert})\nonumber\\
		&+2\rho\sum\limits_{i\in\alpha,\,j\in\gamma}\big(\frac{\lambda_j(A)/\lambda_i(A)}{-\lambda_j(A)/\lambda_i(A)+\rho_3}\big)^2\widetilde{B}_{ij}^2
		-2\sum\limits_{i\in\alpha,\,j\in\gamma}\frac{\lambda_j(A)}{\lambda_i(A)}(\frac{\rho\widetilde{B}_{ij}}{-\lambda_j(A)/\lambda_i(A)+\rho_3})^2+\rho\|b\|^2,
	\end{align*}
	where $\widetilde{B}=P^TBP$ with $P\in{\cal O}^n(A)$. 
\end{lemma}
{\bf Proof. }From \cite[Theorem 6.2]{mohammadi20} and \cite[Lemma 3.1]{Sun06}, we know that 
$${\rm d}^2\delta_{\cal K}(\Phi({x}),\zeta)(z)
=-\varUpsilon_{G({x})}\big({\Gamma},Z\big)+\delta_{{\cal C}_{{\cal K}}}(\Phi({x}),\zeta)(z),$$
where $z=(w,Z)\in\mathcal{Y}\times{\cal S}^n$ and $\varUpsilon_{G({x})}\big({\Gamma},Z\big)=2\langle{\Gamma},ZG({x})^{\dagger}Z\rangle$ is the $\sigma$-term. 
Combining this with the definition of Moreau envelop in \eqref{eqdef:more}, we have 
\begin{align}
	&e_{1/(2\rho)}\big({\rm d}^2\delta_{{\cal K}}(\Phi({x}),\zeta)\big)(b,B)=\inf\limits_{z}\{\rho\|z-(b,B)\|^2+{\rm d}^2\delta_{\cal K}(\Phi({x}),\zeta)(z)\}\nonumber\\
	&=\inf\limits_{z\in{\cal C}_{{\cal K}}(\Phi({x}),{\zeta})}\{-\varUpsilon_{G({x})}\big({\Gamma},Z\big)+\rho\|z-(b,B)\|^2\}\nonumber\\
	&=\inf\limits_{Z\in{\cal C}_{{\cal S}_+^n}(G({x}),{\Gamma})}\{-\varUpsilon_{G({x})}\big({\Gamma},Z\big)+\rho\|Z-B\|^2\}+\rho\|b\|^2. \label{eq:subsetproof9}
\end{align}
From \cite[(28)]{Sun06}, we have
$$-\varUpsilon_{G({x})}\big({\Gamma},Z\big)=-2\sum\limits_{i\in\alpha,\,j\in\gamma}\frac{\lambda_j(A)}{\lambda_i(A)}(P^TZP)_{ij}^2.$$
From \cite[(19)]{Sun06}, we also have
$${\cal C}_{{\cal S}_+^n}(G({x}),{\Gamma})=\left\{Z \in \mathcal{S}^{n}: P_{\beta}^{T} Z P_{\beta} \in{\cal S}_+^{|\beta|}, P_{\beta}^{T}Z P_{\gamma}=0, P_{\gamma}^{T} Z P_{\gamma}=0\right\}.$$
Let $\widetilde{Z}=P^TZP$, $\widetilde{B}=P^TBP$. Thus 
\begin{align}
	&\inf\limits_{Z\in{\cal C}_{{\cal S}_+^n}(G({x}),{\Gamma})}\{-\varUpsilon_{G({x})}\big({\Gamma},{Z}\big)+\rho\|Z-B\|^2\}\nonumber\\
	&=\inf\limits_{Z\in{\cal C}_{{\cal S}_+^n}(G({x}),{\Gamma})}\{-2\sum\limits_{i\in\alpha,\,j\in\gamma}\frac{\lambda_j(A)}{\lambda_i(A)}(P^TZP)_{ij}^2+\rho\|Z-B\|^2\}\nonumber\\
	&=\inf\limits_{Z\in{\cal C}_{{\cal S}_+^n}(G({x}),{\Gamma})}\{-2\sum\limits_{i\in\alpha,\,j\in\gamma}\frac{\lambda_j(A)}{\lambda_i(A)}\widetilde{Z}_{ij}^2+\rho\sum\limits_{i,j}\lvert\widetilde{Z}_{ij}-\widetilde{B}_{ij}\rvert^2\}\label{eq:calex9}. 
\end{align}
For all $i\in\alpha$, $j\in\alpha\cup\beta$ and $i\in\beta$, $j\in\alpha$, to obtain the minmum of \eqref{eq:calex9}, let $\widetilde{Z}_{ij}=\widetilde{B}_{ij}$. For all $i\in\alpha$, $j\in\gamma$ and $i\in\gamma$, $j\in\alpha$, we get the optimal solution $\widetilde{Z}_{ij}=\frac{\rho\widetilde{B}_{ij}}{-\lambda_j(A)/\lambda_i(A)+\rho}$. For all $i,j\in\beta$, $\widetilde{Z}_{\beta\beta}=\Pi_{{\cal S}_+^{\lvert\beta\rvert}}(\widetilde{B}_{\beta\beta})$. Otherwise, $\widetilde{Z}_{ij}=0$. It follows from \eqref{eq:calex9} that 
\begin{align}
	&\inf\limits_{Z\in{\cal C}_{{\cal S}_+^n}(G({x}),{\Gamma})}\{-\varUpsilon_{G({x})}\big({\Gamma},{Z}\big)+\rho\|Z-B\|^2\}\nonumber\\
	=&2\rho\sum\limits_{i\in\beta,\,j\in\gamma}\widetilde{B}_{ij}^2+\rho\sum\limits_{i,\,j\in\gamma}\widetilde{B}_{ij}^2+\rho\,{\rm dist}^2(\widetilde{B}_{\beta\beta},{\cal S}_+^{\lvert\beta\rvert})\nonumber\\
	&+2\rho\sum\limits_{i\in\alpha,\,j\in\gamma}\big(\frac{\lambda_j(A)/\lambda_i(A)}{-\lambda_j(A)/\lambda_i(A)+\rho}\big)^2\widetilde{B}_{ij}^2
	-2\sum\limits_{i\in\alpha,\,j\in\gamma}\frac{\lambda_j(A)}{\lambda_i(A)}(\frac{\rho\widetilde{B}_{ij}}{-\lambda_j(A)/\lambda_i(A)+\rho})^2\label{eq:pr19}
\end{align}
Combining \eqref{eq:subsetproof9} and \eqref{eq:pr19} we have completed the proof. 
$\BOX$
\vskip 10pt
The following result illustrates the uniform second-order expansion for $e\delta_{{\cal S}_+^n}(\cdot)$, which will be used in the derivation of Proposition \ref{prop:uni-soel}. It is worth to note that the second-order expansion for it is firstly studied in \cite[Theorem 3.5]{PRock}. 
By taking advantage of Lemmas \ref{lemma:unih}-\ref{lemma:morcal}, we can provide a direct proof here. 
\begin{proposition}\label{prop:more-souni}
Given $(G(\bar{x}),\overline{\Gamma})\in{\rm gph}\,N_{{\cal S}_+^n}$. Denote $\overline{A}=G(\bar{x})+\overline{\Gamma}$ and $\overline{A}$ possesses the eigenvalue decomposition in \eqref{eq:eig-decomp} with $\overline{P}\in{\cal O}^n(\overline{A})$.	Let $0<r<\min_{i<j}\{v_i(\overline{A})-v_j(\overline{A})\}/3$. For any $A:=G(\bar{x})+{\Gamma}\in\mathbb{B}_r(\overline{A})$ with $\Gamma\in N_{{\cal S}_+^n}(G(\bar{x}))$ and $\pi(A)=\pi(\overline{A})$, we have for all $H\rightarrow0$, 
	$$e\delta_{{\cal S}_+^n}(A+H)-e\delta_{{\cal S}_+^n}(A)=\langle\Pi_{{\cal S}_-^n}(A),H\rangle+\frac{1}{2}e\big({\rm d}^2\delta_{{\cal S}_+^n}(G(\bar{x}),\Gamma)\big)(H)+{O}(\|H\|^3),$$
	where ${O}(\|H\|^3)$ is uniform for all $A\in\mathbb{B}_r(\overline{A})$ with $\Gamma\in N_{{\cal S}_+^n}(G(\bar{x}))$ and $\pi(A)=\pi(\overline{A})$, ${\rm d}^2\delta_{{\cal S}_+^n}(G(\bar{x}),\Gamma)$ is defined in Definition \ref{eqdef:sosd}. 
\end{proposition}
{\bf Proof.} From \cite[Theorem 2.26]{RoWe98}, we have $\nabla e\delta_{{\cal S}_+^n}(A)=\Pi_{{\cal S}_-^n}(A)$. Denote $Q(A)=\Pi_{{\cal S}_-^n}(A)$. 
It follows that 
\begin{align}
	&e\delta_{{\cal S}_+^n}(A+H)-e\delta_{{\cal S}_+^n}(A)=\frac{1}{2}\langle Q(A+H)-Q(A),Q(A+H)+Q(A)\rangle\nonumber\\
	=&\frac{1}{2}\langle H,Q(A)+Q(A+H)-Q(A)\rangle+\frac{1}{2}\langle H,Q(A)\rangle\nonumber\\
	&-\frac{1}{2}\langle\Pi_{{\cal S}_+^n}(A+H)-\Pi_{{\cal S}_+^n}(A),Q(A+H)\rangle-\frac{1}{2}\langle\Pi_{{\cal S}_+^n}(A+H)-\Pi_{{\cal S}_+^n}(A),Q(A)\rangle\nonumber\\
	=&\langle H,Q(A)\rangle+\frac{1}{2}\langle H,Q'(A,H)+{O}(\|H\|^2)\rangle \nonumber\\
	&-\frac{1}{2}\langle\Pi_{{\cal S}_+^n}(A+H)-\Pi_{{\cal S}_+^n}(A),Q(A+H)\rangle-\frac{1}{2}\langle\Pi_{{\cal S}_+^n}(A+H)-\Pi_{{\cal S}_+^n}(A),Q(A)\rangle\label{eq:p5prun}\\
	=&\langle H,Q(A)\rangle+\frac{1}{2}\langle H,Q'(A,H)\rangle+{O}(\|H\|^3)\nonumber\\
	&-\frac{1}{2}\langle\Pi_{{\cal S}_+^n}(A+H)-\Pi_{{\cal S}_+^n}(A),Q(A+H)\rangle-\frac{1}{2}\langle\Pi_{{\cal S}_+^n}(A+H)-\Pi_{{\cal S}_+^n}(A),Q(A)\rangle\label{eq:p5pr0}
\end{align}
where \eqref{eq:p5prun} comes from Proposition \ref{prop:unisocpi}. We denote $\lambda_i:=\lambda_i(A)$, $\Lambda:=\Lambda(A)$ and $\Xi:=\Lambda(A+H)$ for short. Let $P\in{\cal O}^n(A)$, $\widetilde{H}=P^THP$ and $E\in{\cal S}^n$ be the matrix whose components are all 1. 
It follows from $H=\Pi'_{{\cal S}_+^n}(A,H)+\Pi'_{{\cal S}_-^n}(A,H)$, \eqref{eq:dd-projection} and \eqref{eq:dd-projection-} that 
\begin{align}
	&\langle H,Q'(A,H)\rangle=\langle\Pi'_{{\cal S}_+^n}(A,H),\Pi'_{{\cal S}_-^n}(A,H)\rangle+\|\Pi'_{{\cal S}_-^n}(A,H)\|^2\nonumber\\
	&=\langle \Sigma_{\alpha\gamma} \circ {\widetilde H}_{\alpha\gamma}, (E-{\Sigma})_{\alpha\gamma}^T\rangle+\langle \Pi_{{\cal S}_+^{\rvert\beta\lvert} }({\widetilde H}_{\beta\beta}), \Pi_{{\cal S}_-^{\rvert\beta\lvert} }({\widetilde H}_{\beta\beta})\rangle+\langle {\Sigma}_{\alpha\gamma}^T \circ {\widetilde H}_{\alpha\gamma}^T, (E-\Sigma)_{\alpha\gamma} \circ {\widetilde H}_{\alpha\gamma}\rangle+\|\Pi'_{{\cal S}_-^n}(A,H)\|^2\nonumber\\
	&=2\sum_{i\in\alpha,j\in\gamma}\frac{-\lambda_i\lambda_j}{(\lambda_i-\lambda_j)^2}\widetilde{H}_{ij}^2+\|\Pi'_{{\cal S}_-^n}(A,H)\|^2\nonumber\\
	&=2\sum_{i\in\alpha,j\in\gamma}\frac{\lambda_j^2-\lambda_i\lambda_j}{(\lambda_i-\lambda_j)^2}\widetilde{H}_{ij}^2+\sum_{i,j\in\gamma}\widetilde{H}_{ij}^2+2\sum_{i\in\beta,j\in\gamma}\widetilde{H}_{ij}^2+\|\Pi_{{\cal S}_-^{\lvert\beta\rvert}}(\widetilde{H}_{\beta\beta})\|^2\nonumber\\
	&=2\sum_{i\in\beta,j\in\gamma}\widetilde{H}_{ij}^2+\sum_{i,j\in\gamma}\widetilde{H}_{ij}^2+{\rm dist}^2(\widetilde{H}_{\beta\beta},{\cal S}_+^{\lvert\beta\rvert})+2\sum_{i\in\alpha,j\in\gamma}\frac{\lambda_j^2-\lambda_i\lambda_j}{(\lambda_i-\lambda_j)^2}\widetilde{H}_{ij}^2\nonumber\\
	&=\inf\limits_{z}\{\|z-H\|^2+{\rm d}^2\delta_{\cal K}(G(\bar{x}),\Gamma)(z)\}, \label{eq:sigtp}
\end{align}
where the last equality comes from Lemma \ref{lemma:morcal}. 
It is easy to see from the projection onto SDP \eqref{eq:eig-decomp}  that 
\begin{align}\label{eq:p5pr1}
	&-\langle \Pi_{{\cal S}_+^n}(A+H)-\Pi_{{\cal S}_+^n}(A),Q(A+H)\rangle=\langle\Pi_{{\cal S}_+^n}(A),Q(A+H)\rangle\nonumber\\
	&=\langle P\left[\begin{array}{ccc}
		\Lambda(A)_{\alpha\alpha} & 0 & 0 \\ [3pt]
		0 & 0 & 0\\ [3pt]
		0 & 0 & 0
	\end{array}\right]P^T, U\Xi_- U^T\rangle\nonumber\\
	&={\rm tr}(\Lambda_{\alpha\alpha}(P^TU)_{\alpha\gamma}\Xi_{\gamma\gamma}(P^TU)_{\alpha\gamma}^T)+{\rm tr}(\Lambda_{\alpha\alpha}(P^TU)_{\alpha\beta}\Xi_{\beta\beta}(P^TU)_{\alpha\beta}^T),
\end{align}
where $U\in{\cal O}^n(A+H)$ and $(\Xi_-)_{ij}=\min\{0,\Xi_{ij}\}$. 
It follows from Lemma \ref{lemma:unisec} and the Lipschitz continuity of $\lambda(\cdot)$ that 
\begin{align*}
	(P^TU)_{\alpha\gamma}\Xi_{\gamma\gamma}(P^TU)_{\alpha\gamma}^T&=\big[\Theta_{\alpha\gamma}\circ(\widetilde{H}_{\alpha\gamma}Q_{\gamma})\big]\Xi_{\gamma\gamma}\big[\Theta_{\alpha\gamma}^T\circ(Q_{\gamma}^T\widetilde{H}_{\alpha\gamma}^T)\big]+{O}(\|H\|^3)\\
	&=\big[\Theta_{\alpha\gamma}\circ(\widetilde{H}_{\alpha\gamma}Q_{\gamma})\big]\Lambda_{\gamma\gamma}\big[\Theta_{\alpha\gamma}^T\circ(Q_{\gamma}^T\widetilde{H}_{\alpha\gamma}^T)\big]+{O}(\|H\|^3). 
\end{align*}
Thus we have 
$${\rm tr}(\Lambda_{\alpha\alpha}(P^TU)_{\alpha\gamma}\Xi_{\gamma\gamma}(P^TU)_{\alpha\gamma}^T)=\sum_{i\in\alpha,j\in\gamma}\frac{\lambda_i\lambda_j}{(\lambda_i-\lambda_j)^2}\|(\widetilde{H}_{\alpha\gamma}Q_{\gamma})_{ij}\|^2+{O}(\|H\|^3),$$
where ${O}(\|H\|^3)$ is uniform for all $A\in\mathbb{B}_r(\overline{A})$ with $\Gamma\in N_{{\cal S}_+^n}(G(\bar{x}))$ and $\pi(A)=\pi(\overline{A})$. 
Also, it is easy to see that ${\rm tr}(\Lambda_{\alpha\alpha}(P^TU)_{\alpha\beta}\Xi_{\beta\beta}(P^TU)_{\alpha\beta}^T)={O}(\|H\|^3)$. Taking this into \eqref{eq:p5pr1}, we have 
$$-\langle \Pi_{{\cal S}_+^n}(A+H)-\Pi_{{\cal S}_+^n}(A),Q(A+H)\rangle
=\sum_{i\in\alpha,j\in\gamma}\frac{\lambda_i\lambda_j}{(\lambda_i-\lambda_j)^2}\|(\widetilde{H}_{\alpha\gamma}Q_{\gamma})_{ij}\|^2+{O}(\|H\|^3).$$
Similarly, we can compute $\langle\Pi_{{\cal S}_+^n}(A+H)-\Pi_{{\cal S}_+^n}(A),Q(A)\rangle$ in the exactly same way, i.e., 
\begin{align}
	\langle\Pi_{{\cal S}_+^n}(A+H)-\Pi_{{\cal S}_+^n}(A),Q(A)\rangle&={\rm tr}(\Lambda_{\gamma\gamma}(P^TU)_{\gamma\alpha}\Xi_{\alpha\alpha}(P^TU)_{\gamma\alpha}^T)+{\rm tr}(\Lambda_{\gamma\gamma}(P^TU)_{\gamma\beta}\Xi_{\beta\beta}(P^TU)_{\gamma\beta}^T)\nonumber\\
	&=\sum_{i\in\alpha,j\in\gamma}\frac{\lambda_i\lambda_j}{(\lambda_i-\lambda_j)^2}\|(\widetilde{H}_{\gamma\alpha}Q_{\alpha})_{ij}\|^2+{O}(\|H\|^3).\label{eq:anouni}
\end{align}
Combining \eqref{eq:p5pr0}, \eqref{eq:sigtp}, \eqref{eq:p5pr1} and \eqref{eq:anouni} together, we have obtained the result. 
$\BOX$

\section{ALM convergence for NLSDP}\label{sec5}
In this section, we shall study the convergence of the inexact augmented Lagrangian method (ALM) for NLSDP  \eqref{eq:NLSDP} without requiring the uniqueness of Lagrangian multipliers. 
Recalling the augmented Lagrangian function of NLSDP defined in \eqref{eq:deflag}. 
We define the residual function by 
\begin{equation}\label{eq:residual func}
R(x,\zeta)=\|\nabla_xL(x,\zeta)\|+\|\Phi(x)-\Pi_{\cal K}(\Phi(x)+\zeta)\|.
\end{equation}
The augmented Lagrangian method is stated below. 
\begin{algorithm}
\caption{(Augmented Lagrangian method)}\label{algo1}
\begin{algorithmic}[1]
\REQUIRE Let $(x^0,\zeta^0)\in\mathcal{X}\times\mathcal{Y}\times\mathcal{S}^n$, $\rho^0>0$, $\varsigma>1$, $\{\epsilon_k\}_{k\geq0}$ with $\epsilon_k>0$ for all $k$ and $\epsilon_k\rightarrow0$ and set $k:=0$. 
\STATE If $(x^k,\zeta^k)$ satisfies a suitable termination criterion: STOP. 
\STATE Compute $x^{k+1}$ such that 
\begin{equation}\label{eq:subproblemalm}
\|\nabla_x\mathscr{L}(\cdot,\zeta^k,\rho^k)\|\leq\epsilon_k.
\end{equation} 
\STATE Update the vector of multipliers to 
\begin{equation}\label{eq:lambda update}
\zeta^{k+1}:=\rho^k\big[\Phi(x^{k+1})+\frac{\zeta^k}{\rho^k}-\Pi_{\cal K}(\Phi(x^{k+1})+\frac{\zeta^k}{\rho^k})\big].
\end{equation}
\STATE Update $\rho^{k+1}$ by $\rho^{k+1}=\rho^k$ or $\rho^{k+1}=\varsigma\rho^k$ according to certain rules. 
\STATE Set $k\leftarrow k+1$ and go to 1. 
\end{algorithmic}
\end{algorithm}

Before we go any further, it is necessary to introduce the definition of semi-isolated calmness. 
For given perturbation parameters $(a_1,a_2,b)\in {\cal X}\times\mathcal{Y}\times{\cal S}^n$, the corresponding canonical perturbation counterpart of \eqref{eq:NLSDP} is given by  
\begin{equation}\label{eq:NLSDPp1}
	\begin{array}{cl}
		\displaystyle{\min_{x\in {\cal X}}} & f(x) - \langle a_1,x\rangle \\ [3pt]
		{\rm s.t.} & h(x)-a_2=0, \\ [3pt]
		& G(x)-b\in {\cal S}_+^n.
	\end{array} 
\end{equation}
Denote ${S}_{KKT}(a_1,a_2,b)$ the solution set of the KKT optimality condition for problem \eqref{eq:NLSDPp1}, i.e.,
\begin{equation}\label{eq:KKT-NLSDP-p1}
	{S}_{KKT}(a_1,a_2,b)=\left\{(x,y, \Gamma)\in {\cal X}\times\mathcal{Y}\times{\cal S}^n: \left.
	\begin{array}{l}
		\nabla_x L(x,y,\Gamma)-a_1=0,\\ [3pt]
		h(x)-a_2=0, \\ [3pt]
		{\cal S}_+^n\ni (G(x)-b)\perp \Gamma \in {\cal S}_-^n.
	\end{array}
	\right.\right \}. 
\end{equation}
For any KKT pair $(x,y,\Gamma)$ that satisfies the KKT condition with $({a}_1,{a}_2,{b})=(0,0,0)$, we call $x$ the stationary point. For a given feasible solution $\bar{x}\in {\cal X}$ of \eqref{eq:NLSDP} with $({a}_1,{a}_2,{b})=(0,0,0)$, define ${\cal M}(\bar{x})$ the set of all multiples $(y,\Gamma)\in\mathcal{Y}\times{\cal S}^n$ satisfying the KKT condition \eqref{eq:KKT-NLSDP-p1}, i.e.,
\begin{equation}\label{eq:defcalm}
	{\cal M}(\bar{x})=\{(y,\Gamma)\in{\cal Y}\times{\cal S}^n \mid (\bar{x},y,\Gamma)\in {S}_{\rm KKT}(0,0,0)\}.
\end{equation}
 It is easy to see ${\cal M}(\bar{x})$ is convex. 

The definition of semi-isolated calmness of a set-valued mapping is first officially presented by \cite[Theorem 4.1]{MSarabi2018}, which is an extension of \cite[Proposition 1.43]{ISolodov14}, to establish the characterization of noncritical multipliers for the polyhedral problem (to be more specifically, composite piecewise linear quadratic problem). It is worth noting that for the polyhedral case, semi-isolated calmness is equivalent to noncritical \cite[Theorem 4.1]{MSarabi2018}, though this does not hold for the non-polyhedral case as stated in \cite[Theorem 5.6]{MSarabi19}. 
\begin{definition}\label{def:semi-ic}
	The semi-isolated clamness for the mapping $S_{\rm KKT}$ at $((0,0,0),(\bar{x},\bar{y},\overline{\Gamma}))$ holds if there exists $\kappa>0$ and open neighborhoods $\mathbb{U}$ of $(0,0,0)$ and $\mathbb{V}$ of $(\bar{x},\bar{y},\overline{\Gamma})$ such that for all $(a_1,a_2,b)\in\mathbb{U}$, 
	\begin{equation*}
		\|x-\bar{x}\|+{\rm dist}((y,\Gamma),{\cal M}(\bar{x}))\le \kappa\|(a_1,a_2,b)\|\quad \forall\, (x,y,\Gamma)\in{S}_{\rm KKT}(a_1,a_2,b)\cap \mathbb{V}.
	\end{equation*}
\end{definition}

Let $\bar{x}$ be a stationary point of NLSDP \eqref{eq:NLSDP}. Assume ${\cal M}(\bar{x})$ is nonempty with $(\bar{y},\overline{\Gamma})\in{\cal M}(\bar{x})$, the critical cone of problem \eqref{eq:NLSDP} is adopted from \cite[(37)]{Sun06} 
\begin{align*}
	{\cal C}(\bar{x})
	=&\left\{d\in{\cal X}: \nabla h(\bar{x}) d=0, \nabla G(\bar{x}) d \in {\cal C}_{\mathcal{S}_{+}^{n}}(G(\bar{x}),\overline{\Gamma})\right\}, 
\end{align*}
where ${\cal C}_{\mathcal{S}_{+}^{n}}(G(\bar{x}),\overline{\Gamma})$ is the critical cone defined in \eqref{eq:crisdp}. 
We also need the definition of second order sufficient condition for \eqref{eq:NLSDP}, which can be found from, e.g., \cite[equation (2.11)]{HSarabi20}. 
\begin{definition}\label{def:sosc}
	Let $\bar{x}$ be a stationary point of NLSDP \eqref{eq:NLSDP}. Given $(\bar{y},\overline{\Gamma})\in{\cal M}(\bar{x})$. We say the second order sufficient condition (SOSC) holds at $(\bar{x},\bar{y},\overline{\Gamma})$ if 
	\begin{equation}\label{eq:soscr}
	\langle \nabla_{xx}^2L(\bar{x},\bar{y},\overline{\Gamma})d,d\rangle-\varUpsilon_{G(\bar{x})}\big(\overline{\Gamma},\nabla G(\bar{x}){d}\big)>0,\quad\forall\; 0\neq d\in{\cal C}(\bar{x}).
	\end{equation}
	where $\varUpsilon_{G(\bar{x})}\big(\overline{\Gamma},\nabla G(\bar{x}){d}\big)=2\langle\overline{\Gamma},(\nabla G(\bar{x}){d})G(\bar{x})^{\dagger}(\nabla G(\bar{x}){d})\rangle$ is the $\sigma$-term (cf. \cite[Lemma 3.1]{Sun06}) and $G(\bar{x})^{\dagger}$ is the generalized inverse matrix of $G(\bar{x})$. 
\end{definition}

\subsection{Properties of augmented Lagrangian function and the solution of subproblem \eqref{eq:subproblemalm}}
To establish the ALM convergence of NLSDP, we firstly need the quadratic growth condition for augmented Lagrangian function as presented in \cite{SSZhang08} and \cite{HSarabi20}. The (uniform) positive definite condition is critical in establishing the (uniform) quadratic growth condition. In \cite[Lemma 4.2]{HSarabi20}, they studied the uniform version for linearly quadratic composite optimization problems under only SOSC. We try to extend their result from the polyhedral problem to NLSDP. It is worth noting that similar result for NLSDP was established in \cite[Proposition 4]{SSZhang08}, although they supposed nondegeneracy and strongly SOSC. 
It follows from \cite[Theorem 8.3]{mohammadi20} that for any given $\zeta\in{\cal M}(\bar{x})$, function $x\mapsto\mathscr{L}({x},\zeta,\rho)$ is twice semidifferentiable with respect to $x$ at $\bar{x}$. Then we obtain the following lemma. 
\begin{lemma}\label{lemma:soscegeq0}
	Let $\bar{x}\in{\cal X}$ be a stationary point to the NLSDP \eqref{eq:NLSDP} 
	and $\bar{\zeta}\in {\cal M}(\bar{x})$ \eqref{eq:defcalm}. Then the following conditions are equivalent:
	\begin{itemize}
		\item[(a)] the SOSC holds at $(\bar{x},\bar{\zeta})$ (see Definition \ref{def:sosc}); 
		\item[(b)] there are positive constants $\rho_3$, $\varepsilon'$, $l'$ such that for all $\zeta\in{\cal M}(\bar{x})\cap\mathbb{B}_{\varepsilon'}(\bar{\zeta})$ and all $\rho\geq\rho_3$, 
		\begin{equation}\label{eq:sosceqine}
			{\rm d}_x^2\mathscr{L}\big(\bar{x},\zeta,\rho\big)(w)\geq l'\|w\|^2\quad\forall\, w\in\Re^n\backslash\{0\}, 
		\end{equation}
		where 
		${\rm d}_x^2\mathscr{L}$ is the second semiderivative defined by Definition \ref{eqdef:sosd}. 
	\end{itemize}
\end{lemma}
{\bf Proof.}  ``(b)$\Rightarrow$(a)": It follows from \cite[Theorem 8.4]{mohammadi20}, directly. 

``(a)$\Rightarrow$(b)": First, it follows from \cite[Proposition 13.5]{RoWe98} that the second semiderivative is lower semicontinuous and positive homogenous of degree 2. Thus, by  \cite[Theorem 8.4 (ii)]{mohammadi20}, we know that the SOSC is equivalent to the existence of $l'>0$ such that  
\begin{equation}\label{eq:nonunif}
	{\rm d}_x^2\mathscr{L}\big(\bar{x},\bar{\zeta},\rho_3\big)(w)\geq2l'\quad \forall\,w\in\mathbb{S},
\end{equation}
where $\mathbb{S}$ is the unite sphere of ${\cal X}$. 

Next, we shall show that there exists an $\varepsilon'>0$ such that for any $\zeta\in{\cal M}(\bar{x})\cap\mathbb{ B}_{\varepsilon'}(\bar{\lambda})$, 
$${\rm d}_x^2\mathscr{L}\big(\bar{x},\zeta,\rho_3\big)(w)\geq l'\quad \forall\,w\in\mathbb{S}.$$
From \cite[Theorem 8.3 (i)]{mohammadi20}, we know that for any $\zeta\in{\cal M}(\bar{x})$, 
\begin{align}
	{\rm d}_x^2\mathscr{L}(\bar{x},\zeta,\rho_3)(w)=&\langle\nabla_{xx}^2L(\bar{x},\zeta)w,w\rangle\nonumber\\
	&+\inf\limits_{z}\{\rho_3\|z-\nabla\Phi(\bar{x})w\|^2+{\rm d}^2\delta_{\cal K}(\Phi(\bar{x}),\zeta)(z)\}, \label{eq:profd2}
\end{align}
where  ${\cal K}=\{0\}\times{\cal S}^n_+$. 
Choose $\varepsilon''\in(0,l'/(2\|\nabla^2\Phi(\bar{x})\|))$ if $\|\nabla^2\Phi(\bar{x})\|\neq0$ and $\varepsilon''>0$ otherwise. For each $w\in\mathbb{S}$, it is clear that for any $\zeta\in\mathbb{B}_{\varepsilon''}(\bar{\zeta})$, 
\begin{align}
	\langle\nabla_{xx}^2L(\bar{x},\zeta)w,w\rangle&=\langle\nabla_{xx}^2L(\bar{x},\bar{\zeta})w,w\rangle+\langle\zeta-\bar{\zeta},\nabla^2\Phi(\bar{x})(w,w)\rangle\nonumber\\
	&\geq\langle\nabla_{xx}^2L(\bar{x},\bar{\zeta})w,w\rangle-\|\nabla^2\Phi(\bar{x})\|\cdot\|\zeta-\bar{\zeta}\|\nonumber\\
	&\geq\langle\nabla_{xx}^2L(\bar{x},\bar{\zeta})w,w\rangle-\frac{l'}{2}. \label{eq:lagg}
\end{align}
Let $A=\Gamma+G(\bar{x})$ and $\overline{A}=\overline{\Gamma}+G(\bar{x})$. Suppose $\overline{P}\in{\cal O}^n(\overline{A})$ and $P\in{\cal O}^n(A)$. Suppose $\overline{A}$ possesses the eigenvalue decomposition \eqref{eq:eig-decomp}. Let $0<r<\min_{i<j}\{v_i(\overline{A})-v_j(\overline{A})\}/3$. For all $\Gamma\in\mathbb{B}_r(\overline{\Gamma})$, denote $\alpha\cup\beta_+:=\{i\mid\lambda(A)>0\}$, $\beta_0:=\{i\mid\lambda(A)=0\}$ and $\gamma\cup\beta_-:=\{i\mid\lambda(A)<0\}$. {Since $\Gamma\in{\cal N}_{{\cal S}_+^n}(G(\bar{x}))\cap\mathbb{B}_r(\overline{\Gamma})$, we know that $\Gamma$ possess the following structure:
$$\Gamma=P\left[\begin{array}{ccccc}
		0_{\alpha\alpha} & 0 & 0 & 0 & 0 \\ [3pt]
		0 & 0_{\beta_+\beta_+} & 0 & 0 & 0\\ [3pt]
		0 & 0 & 0_{\beta_0\beta_0} & 0 & 0\\ [3pt]
		0 & 0 & 0 & \Lambda(\Gamma)_{\beta_-\beta_-} & 0\\ [3pt]
		0 & 0 & 0 & 0 & \Lambda(\Gamma)_{\gamma\gamma}
	\end{array}\right]P^T, $$ 
where $\Lambda(\Gamma)_{\beta_-\beta_-}<0$ and $\Lambda(\Gamma)_{\gamma\gamma}<0$. Thus the index set of positive eigenvalues of $A$ remains unchanged as the positive eigenvalues part of $A$ is provided by $G(\bar{x})$, which is unchanged. Then we have $\beta_+=\varnothing$  and $\beta=\beta_0\cup\beta_-$. } Let $\widetilde{Z}=P^TZP$, $\widetilde{B}=P^T(\nabla G(\bar{x})w)P$. It follows from Lemma \ref{lemma:morcal} that  
\begin{align}
	&\inf\limits_{z}\{\rho_3\|z-\nabla\Phi(\bar{x})w\|^2+{\rm d}^2\delta_{\cal K}(\Phi(\bar{x}),\zeta)(z)\}\nonumber\\
	=&\inf\limits_{Z\in{\cal C}_{{\cal S}_+^n}(G(\bar{x}),{\Gamma})}\{-\varUpsilon_{G(\bar{x})}\big({\Gamma},{Z}\big)+\rho_3\|Z-\nabla G(\bar{x})w\|^2\}+\rho_3\|\nabla h(\bar{x})w\|^2\nonumber\\
	=&2\rho_3\sum\limits_{i\in\beta_0,\,j\in\gamma\cup\beta_-}\widetilde{B}_{ij}^2+\rho_3\sum\limits_{i,\,j\in\gamma\cup\beta_-}\widetilde{B}_{ij}^2+\rho_3{\rm dist}^2(\widetilde{B}_{\beta_0\beta_0},{\cal S}_+^{\lvert\beta_0\rvert})\nonumber\\
	&+2\rho_3\sum\limits_{i\in\alpha,\,j\in\gamma\cup\beta_-}\big(\frac{\lambda_j(A)/\lambda_i(A)}{-\lambda_j(A)/\lambda_i(A)+\rho_3}\big)^2\widetilde{B}_{ij}^2\nonumber\\
	&-2\sum\limits_{i\in\alpha,\,j\in\gamma\cup\beta_-}\frac{\lambda_j(A)}{\lambda_i(A)}(\frac{\rho_3\widetilde{B}_{ij}}{-\lambda_j(A)/\lambda_i(A)+\rho_3})^2+\rho_3\|\nabla h(\bar{x})w\|^2. \label{eq:pr1}
\end{align}
From \cite[Proposition 2.6]{Ding14}, we know that for $\Gamma$ sufficiently close to $\overline{\Gamma}$, we have 
${\rm dist}(P,{\cal O}^n(\overline{A}))=O(\|\Gamma-\overline{\Gamma}\|),$
which implies for every $A$, there exists $Q={\rm Diag}(Q_1,\dots,Q_{\bar{d}})$,  such that 
\begin{equation}\label{eq:conveigvec}
P=\overline{P}Q+O(\|\Gamma-\overline{\Gamma}\|),
\end{equation}
where $\bar{d}$ is the number of the different eigenvalues of $\overline{A}$ and $Q_l\in{\cal O}^{\lvert\bar{\iota}_l\rvert}$. 
We can take $\Gamma$ such that $\|\Gamma-\overline{\Gamma}\|\leq r$. Let $\widehat{B}=\overline{P}^T(\nabla G(\bar{x})w)\overline{P}$ and ${\cal C}_{{\cal S}_+^{\lvert\beta\rvert}}:=\{W\in{\cal S}_+^{\lvert\beta\rvert}\mid \overline{P}^T_{\beta_0}W\overline{P}_{\beta_0}\in{\cal S}_+^{\lvert\beta\rvert},\;\overline{P}^T_{\beta_0}W\overline{P}_{\beta_-}=0,\;\overline{P}^T_{\beta_-}W\overline{P}_{\beta_-}=0\}$. 
Then we have  
\begin{align*}
	&2\rho_3\sum\limits_{i\in\beta_0,\,j\in\gamma\cup\beta_-}\widetilde{B}_{ij}^2+\rho_3\sum\limits_{i,\,j\in\gamma\cup\beta_-}\widetilde{B}_{ij}^2+\rho_3{\rm dist}^2(\widetilde{B}_{\beta_0\beta_0},{\cal S}_+^{\lvert\beta_0\rvert})\nonumber\\
	=&2\rho_3\sum\limits_{i\in\beta,\,j\in\gamma}\widetilde{B}_{ij}^2+\rho_3\sum\limits_{i,\,j\in\gamma}\widetilde{B}_{ij}^2+\rho_3{\rm dist}^2(\widetilde{B}_{\beta_0\beta_0},{\cal S}_+^{\lvert\beta_0\rvert})+2\rho_3\sum\limits_{i\in\beta_0,\,j\in\beta_-}\widetilde{B}_{ij}^2+\rho_3\sum\limits_{i,\,j\in\beta_-}\widetilde{B}_{ij}^2\nonumber\\
	\overset{\eqref{eq:crisdp}}{=}&2\rho_3\sum\limits_{i\in\beta,\,j\in\gamma}\widetilde{B}_{ij}^2+\rho_3\sum\limits_{i,\,j\in\gamma}\widetilde{B}_{ij}^2+\rho_3{\rm dist}^2\big(\widetilde{B}_{\beta\beta},,{\cal C}_{{\cal S}_+^{\lvert\beta\rvert}}\big)
	=2\rho_3\|\widetilde{B}_{\beta\gamma}\|^2+\rho_3\|\widetilde{B}_{\gamma\gamma}\|^2+\rho_3{\rm dist}^2\big(\widetilde{B}_{\beta\beta},{\cal C}_{{\cal S}_+^{\lvert\beta\rvert}}\big)\nonumber\\
	=&2\rho_3\|\widehat{B}_{\beta\gamma}\|^2+\rho_3\|\widehat{B}_{\gamma\gamma}\|^2+\rho_3{\rm dist}^2\big(\widetilde{B}_{\beta\beta},{\cal C}_{{\cal S}_+^{\lvert\beta\rvert}}\big)+O(\|\Gamma-\overline{\Gamma}\|), 
	\end{align*}
where the last equality comes from the diagonal structure of $Q$ and  \eqref{eq:conveigvec}. It follows that 
	\begin{align}
	&2\rho_3\sum\limits_{i\in\beta_0,\,j\in\gamma\cup\beta_-}\widetilde{B}_{ij}^2+\rho_3\sum\limits_{i,\,j\in\gamma\cup\beta_-}\widetilde{B}_{ij}^2+\rho_3{\rm dist}^2(\widetilde{B}_{\beta_0\beta_0},{\cal S}_+^{\lvert\beta_0\rvert})\nonumber\\
	=&2\rho_3\sum\limits_{i\in\beta,\,j\in\gamma}\widehat{B}_{ij}^2+\rho_3\sum\limits_{i,\,j\in\gamma}\widehat{B}_{ij}^2+\rho_3{\rm dist}^2(\widetilde{B}_{\beta\beta},{\cal C}_{{\cal S}_+^{\lvert\beta\rvert}})+O(\|\Gamma-\overline{\Gamma}\|)\nonumber\\	
	\geq&\rho_3{\rm dist}^2\big(\widetilde{B}_{\beta\beta},{\cal S}_+^{\lvert\beta\rvert}\big)+2\rho_3\sum\limits_{i\in\beta,\,j\in\gamma}\widehat{B}_{ij}^2+\rho_3\sum\limits_{i\in\gamma,\,j\in\gamma}\widehat{B}_{ij}^2+O(\|\Gamma-\overline{\Gamma}\|)\nonumber\\
	=&\rho_3{\rm dist}^2\big(\widehat{B}_{\beta\beta},{\cal S}_+^{\lvert\beta\rvert}\big)+2\rho_3\sum\limits_{i\in\beta,\,j\in\gamma}\widehat{B}_{ij}^2+\rho_3\sum\limits_{i\in\gamma,\,j\in\gamma}\widehat{B}_{ij}^2+O(\|\Gamma-\overline{\Gamma}\|),\label{eq:gammabar}
\end{align}
where the inequality comes from ${\cal S}_+^{\lvert\beta\rvert}\supseteq{\cal C}_{{\cal S}_+^{\lvert\beta\rvert}}$. 
and the last equality comes from ${\rm dist}\big(\widetilde{B}_{\beta\beta},{\cal S}_+^{\lvert\beta\rvert}\big)=\|\widetilde{B}_{\beta\beta}-\Pi_{{\cal S}_+^{\lvert\beta\rvert}}(\widetilde{B}_{\beta\beta})\|=\|\widehat{B}_{\beta\beta}-\Pi_{{\cal S}_+^{\lvert\beta\rvert}}(\widehat{B}_{\beta\beta})\|+O(\|\Gamma-\overline{\Gamma}\|)={\rm dist}\big(\widehat{B}_{\beta\beta},{\cal S}_+^{\lvert\beta\rvert}\big)+O(\|\Gamma-\overline{\Gamma}\|)$. 

We also have
		\begin{align*}
			&2\rho_3\sum\limits_{i\in\alpha,\,j\in\gamma\cup\beta_-}\big(\frac{\lambda_j(A)/\lambda_i(A)}{-\lambda_j(A)/\lambda_i(A)+\rho_3}\big)^2\widetilde{B}_{ij}^2-2\sum\limits_{i\in\alpha,\,j\in\gamma\cup\beta_-}\frac{\lambda_j(A)}{\lambda_i(A)}\big(\frac{\rho_3\widetilde{B}_{ij}}{-\lambda_j(A)/\lambda_i(A)+\rho_3}\big)^2\nonumber\\
			\geq&2\rho_3\sum\limits_{i\in\alpha,\,j\in\gamma}\big(\frac{\lambda_j(\Gamma)}{-\lambda_j(\Gamma)+\rho_3\lambda_i(G(\bar{x}))}\big)^2\widetilde{B}_{ij}^2-2\sum\limits_{i\in\alpha,\,j\in\gamma}\frac{\lambda_j(\Gamma)}{\lambda_i(G(\bar{x}))}\big(\frac{\rho_3\widetilde{B}_{ij}\lambda_i(G(\bar{x}))}{-\lambda_j(\Gamma)+\rho_3\lambda_i(G(\bar{x}))}\big)^2\nonumber\\
			=&2\rho_3\sum\limits_{i\in\alpha,\,j\in\gamma}[\big(\frac{\lambda_j(\overline{\Gamma})}{-\lambda_j(\overline{\Gamma})+\rho_3\lambda_i(G(\bar{x}))})^2+O(\|\Gamma-\overline{\Gamma}\|)](\widetilde{B}_{ij})^2\nonumber\\
			&-2\sum\limits_{i\in\alpha,\,j\in\gamma}[\frac{\rho_3^2\lambda_i(G(\bar{x}))\lambda_j(\overline{\Gamma})}{(-\lambda_j(\overline{\Gamma})+\rho_3\lambda_i(G(\bar{x})))^2}+O(\|\Gamma-\overline{\Gamma}\|)](\widetilde{B}_{ij})^2\nonumber\\
			=&2\rho_3\sum\limits_{i\in\alpha,\,j\in\gamma}[\frac{\lambda_j(\overline{\Gamma})^2-\rho_3\lambda_j(\overline{\Gamma})\lambda_i(G(\bar{x}))}{(-\lambda_j(\overline{\Gamma})+\rho_3\lambda_i(G(\bar{x})))^2}+O(\|\Gamma-\overline{\Gamma}\|)](\widetilde{B}_{ij})^2\nonumber\\
\end{align*}
where the first equality follows from the Lipschitz continuity of $\lambda(\cdot)$ and the locally Lipschitz continuity of $g_1(z)=\big(\frac{z}{-z+\rho_3\lambda_i(G(\bar{x}))})^2$, $g_2(z)=\frac{\rho_3^2\lambda_i(G(\bar{x}))z}{(-z+\rho_3\lambda_i(G(\bar{x})))^2}$ over $z<0$. 
Suppose $\overline{A}$ has $\bar{d}$ different eigenvalues with  $v_1(\overline{A})>\dots>v_{d_0}(\overline{A})>0=v_{d_0+1}(\overline{A})>v_{d_0+2}(\overline{A})>v_{\bar{d}}(\overline{A})$ and denote $\iota_s=\{i\mid\lambda_i(\overline{A})=v_s(\overline{A})\}$. Then we have 
\begin{align}
&2\rho_3\sum\limits_{i\in\alpha,\,j\in\gamma}[\frac{\lambda_j(\overline{\Gamma})^2-\rho_3\lambda_j(\overline{\Gamma})\lambda_i(G(\bar{x}))}{(-\lambda_j(\overline{\Gamma})+\rho_3\lambda_i(G(\bar{x})))^2}+O(\|\Gamma-\overline{\Gamma}\|)](\widetilde{B}_{ij})^2\nonumber\\
&=2\rho_3\sum\limits_{s=1}^{d_0}\sum\limits_{t=d_0+2}^{\bar{d}}[\frac{v_t(\overline{\Gamma})^2-\rho_3v_t(\overline{\Gamma})v_s(G(\bar{x}))}{(-v_t(\overline{\Gamma})+\rho_3v_s(G(\bar{x})))^2}+O(\|\Gamma-\overline{\Gamma}\|)]\sum\limits_{i\in\iota_s,j\in\iota_t}(\widetilde{B}_{ij})^2\nonumber\\
&=2\rho_3\sum\limits_{s=1}^{d_0}\sum\limits_{t=d_0+2}^{\bar{d}}[\frac{v_t(\overline{\Gamma})^2-\rho_3v_t(\overline{\Gamma})v_s(G(\bar{x}))}{(-v_t(\overline{\Gamma})+\rho_3v_s(G(\bar{x})))^2}]\sum\limits_{i\in\iota_s,j\in\iota_t}(\widehat{B}_{ij})^2+O(\|\Gamma-\overline{\Gamma}\|)\nonumber\\ 
&=2\rho_3\sum\limits_{i\in\alpha,\,j\in\gamma}[\frac{\lambda_j(\overline{\Gamma})^2-\rho_3\lambda_j(\overline{\Gamma})\lambda_i(G(\bar{x}))}{(-\lambda_j(\overline{\Gamma})+\rho_3\lambda_i(G(\bar{x})))^2}](\widehat{B}_{ij})^2+O(\|\Gamma-\overline{\Gamma}\|), \label{eq:lambpert}
\end{align}
where the last second inequality follows from the block diagonal structure of $Q$ as $\sum\limits_{i\in\iota_s,j\in\iota_t}(\widetilde{B}_{ij})^2=\|Q_{s}\widehat{B}_{\iota_s\iota_t}Q_{t}^T\|^2=\|\widehat{B}_{\iota_s\iota_t}\|^2=\sum\limits_{i\in\iota_s,j\in\iota_t}(\widehat{B}_{ij})^2$.

It follows from Lemma \ref{lemma:morcal} again that if we combine \eqref{eq:gammabar} and \eqref{eq:lambpert} together, the right hand side is exactly 
$$\inf\limits_{Z\in{\cal C}_{{\cal S}_+^n}(G(\bar{x}),\overline{\Gamma})}\{-\varUpsilon_{G(\bar{x})}\big(\overline{\Gamma},{Z}\big)+\rho_3\|Z-\nabla G(\bar{x})w\|^2\}+\Delta(\Gamma),$$
where $\Delta(\Gamma)=O(\|\Gamma-\overline{\Gamma}\|)$ is number, which satisfies that there exist $q>0$ and $\omega>0$ such that for all $\Gamma\in\mathbb{B}_{\omega}(\overline{\Gamma})$, 
$\lvert \Delta(\Gamma)\rvert\leq q\|\Gamma-\overline{\Gamma}\|$ as it originates from \eqref{eq:conveigvec}. 
Combining this fact with \eqref{eq:pr1}, we have  for all $\Gamma$ such that $\|\Gamma-\overline{\Gamma}\|\leq\min\{r,l'/q,\omega\}$, 
\begin{align}
	&\inf\limits_{z}\{\rho_3\|z-\nabla\Phi(\bar{x})w\|^2+{\rm d}^2\delta_{\cal K}(\Phi(\bar{x}),\zeta)(z)\}\nonumber\\
	&=\inf\limits_{Z\in{\cal C}_{{\cal S}_+^n}(G(\bar{x}),{\Gamma})}\{-\varUpsilon_{G(\bar{x})}\big({\Gamma},{Z}\big)+\rho_3\|Z-\nabla G(\bar{x})w\|^2\}+\rho_3\|\nabla h(\bar{x})w\|^2\nonumber\\
	&\geq\inf\limits_{Z\in{\cal C}_{{\cal S}_+^n}(G(\bar{x}),\overline{\Gamma})}\{-\varUpsilon_{G(\bar{x})}\big(\overline{\Gamma},{Z}\big)+\rho_3\|Z-\nabla G(\bar{x})w\|^2\}+\rho_3\|\nabla h(\bar{x})w\|^2+O(\|\Gamma-\overline{\Gamma}\|)\nonumber\\
	&\geq\inf\limits_{z}\{\rho_3\|z-\nabla\Phi(\bar{x})w\|^2+{\rm d}^2\delta_{\cal K}(\Phi(\bar{x}),\bar{\zeta})(z)\}-l'.\label{eq:myine}
\end{align}
Let $\varepsilon'=\min\{\varepsilon'',l'/q, r,\omega\}$. By \eqref{eq:nonunif}, \eqref{eq:profd2}, \eqref{eq:lagg} and \eqref{eq:myine}, we have verified that for any $\zeta\in{\cal M}(\bar{x})\cap\mathbb{ B}_{\varepsilon'}(\bar{\lambda})$, 
$${\rm d}_x^2\mathscr{L}\big(\bar{x},\zeta,\rho_3\big)(w)\geq{\rm d}_x^2\mathscr{L}\big(\bar{x},\bar{\zeta},\rho_3\big)(w)-l'\geq l' \quad \forall\,w\in\mathbb{S}.$$

Using the positive homogeneity of the second semiderivative yields \eqref{eq:sosceqine} for $\rho=\rho_3$ and for all $\zeta\in{\cal M}(\bar{x})\cap{\cal B}_{\varepsilon'}(\bar{\zeta})$. Recall that the function
$$\rho\mapsto e_{1/2\rho}\big({\rm d}^2\delta_{{\cal K}}(\Phi(\bar{x}),\zeta)\big)(\nabla \Phi(\bar{x})w)$$
is nondecreasing. Therefor the function $\rho\mapsto{\rm d}_x^2\mathscr{L}\big((\bar{x},\bar{\zeta},\rho),0\big)(w)$ is also nondecreasing. This yields \eqref{eq:sosceqine} for all $\zeta\in{\cal M}(\bar{x})\cap{\cal B}_{\varepsilon'}(\bar{\zeta})$ and all $\rho\geq\rho_3$, and hence complete the proof. 
$\BOX$

We also need the following result on the uniform expansion of augmented Lagrangian function, which can be obtained from Proposition \ref{prop:more-souni}. 
\begin{proposition}\label{prop:uni-soel}
	Let $\bar{x}\in{\cal X}$ be a stationary point to the NLSDP \eqref{eq:NLSDP} and $\bar{\zeta}\in{\cal M}(\bar{x})$. Let $\overline{A}=G(\bar{x})+\overline{\Gamma}$ and $0<r<\min_{i<j}\{v_i(\overline{A})-v_j(\overline{A})\}/3$. For all $\zeta\in{\cal M}(\bar{x})\cap\mathbb{B}_{r}(\bar{\lambda})$ with $\pi(\Gamma)=\pi(\overline{\Gamma})$ and any $\rho>0$, we have
	\begin{equation}\label{eq:soel}
		\frac{f(\bar{x})-\mathscr{L}(x,\zeta,\rho)}{\|x-\bar{x}\|^2}=-\frac{1}{2}{\rm d}_x^2\mathscr{L}(\bar{x},\zeta,\rho)(\frac{x-\bar{x}}{\|x-\bar{x}\|})+{O}(\|x-\bar{x}\|),
	\end{equation}
	where ${O}(\|x-\bar{x}\|)$ is uniform for all $\zeta\in{\cal M}(\bar{x})\cap\mathbb{B}_{r}(\bar{\lambda})$ with $\pi(\Gamma)=\pi(\overline{\Gamma})$. 
\end{proposition}
{\bf Proof.} From \cite[Proposition 3.2]{HSarabi20}, we know that for all $\zeta\in{\cal M}(\bar{x})$ and any $\rho>0$, $f(\bar{x})=\mathscr{L}(\bar{x},\zeta,\rho)$. It follows that 
\begin{align}
	\mathscr{L}(x,\zeta,\rho)-f(\bar{x})=&f(x)-f(\bar{x})+\langle y,h(x)\rangle+\frac{\rho}{2}\|h(x)\|^2-\big(\langle y,h(\bar{x})\rangle+\frac{\rho}{2}\|h(\bar{x})\|^2\big)\nonumber\\
	&+\rho[e\delta_{{\cal S}_+^n}(G({x})+\rho^{-1}\Gamma)-e\delta_{{\cal S}_+^n}(G(\bar{x})+\rho^{-1}\Gamma)].\label{eq:diffl1}
\end{align}
It can be checked directly from Proposition \ref{prop:more-souni} that
\begin{align}
	&e\delta_{{\cal S}_+^n}(G({x})+\rho^{-1}\Gamma)-e\delta_{{\cal S}_+^n}(G(\bar{x})+\rho^{-1}\Gamma)-\langle\nabla_xe\delta_{{\cal S}_+^n}(G(\bar{x})+\rho^{-1}\Gamma) ,x-\bar{x}\rangle\nonumber\\
	=&\langle\Pi_{{\cal S}_-^n}(G(\bar{x})+\rho^{-1}\Gamma),G(x)-G(\bar{x})\rangle+\frac{1}{2}e\big({\rm d}^2\delta_{{\cal S}_+^n}(G(\bar{x}),\Gamma)\big)\big(G(x)-G(\bar{x})\big)\nonumber\\
	&-\langle\nabla(e\delta_{{\cal S}_+^n})(G(\bar{x})+\rho^{-1}\Gamma),\nabla G(\bar{x})(x-\bar{x})\rangle+{O}(\|G(x)-G(\bar{x})\|^3) \nonumber \\
	=&\langle\Pi_{{\cal S}_-^n}(G(\bar{x})+\rho^{-1}\Gamma),\frac{1}{2}\nabla^2G(\bar{x})(x-\bar{x},x-\bar{x})\rangle+{O}(\|G(x)-G(\bar{x})\|^3)\nonumber\\
	&+{O}(\|x-\bar{x}\|^3)+\frac{1}{2}e\big({\rm d}^2\delta_{{\cal S}_+^n}(G(\bar{x}),\Gamma)\big)\big(\nabla G(\bar{x})(x-\bar{x})+{O}(\|x-\bar{x}\|^2)\big). \label{eq:diffenv}
\end{align} 
From the explicit form of $e\big({\rm d}^2\delta_{{\cal S}_+^n}(G(\bar{x}),\Gamma)\big)(\cdot)$ in Lemma \ref{lemma:morcal}, we know that 
\begin{align}
	&e\big({\rm d}^2\delta_{{\cal S}_+^n}(G(\bar{x}),\Gamma)\big)\big(\nabla G(\bar{x})(x-\bar{x})+{O}(\|x-\bar{x}\|^2)\big)\nonumber\\ 
	&=e\big({\rm d}^2\delta_{{\cal S}_+^n}(G(\bar{x}),\Gamma)\big)\big(\nabla G(\bar{x})(x-\bar{x})\big)+{O}(\|x-\bar{x}\|^3),\label{eq:diffl}
\end{align}
where ${O}(\|x-\bar{x}\|^2)$ and ${O}(\|x-\bar{x}\|^3)$ in \eqref{eq:diffenv} and \eqref{eq:diffl} are uniform for all $\zeta\in{\cal M}(\bar{x})\cap\mathbb{B}_{r}(\bar{\lambda})$ with $\pi(\Gamma)=\pi(\overline{\Gamma})$. 
Combining the continuity of $f(x)-f(\bar{x})+\langle y,h(x)\rangle+\frac{\rho}{2}\|h(x)\|^2-\big(\langle y,h(\bar{x})\rangle+\frac{\rho}{2}\|h(\bar{x})\|^2\big)$ on $x$ with \eqref{eq:diffl1}, \eqref{eq:diffenv}, \eqref{eq:diffl} and \cite[Theorem 8.3]{mohammadi20} with $\bar{\zeta}\in{\cal M}(\bar{x})$, we have attained \eqref{eq:soel}.
$\BOX$

Combining Lemma \ref{lemma:soscegeq0} and Proposition \ref{prop:uni-soel} together, we are ready to state the uniform quadratic growth condition for augmented Lagrangian function under SOSC. The non-uniform form for NLSDP is firstly studied in \cite[Proposition 1]{SSZhang08} and extended to general $C^2$-reducible constrained optimization by Mohammadi et al. \cite[Theorem 8.4]{mohammadi20} under weaker condition. 
\begin{theorem}\label{thm:equuqgc}
	Let $\bar{x}\in{\cal X}$ be a stationary point to the NLSDP \eqref{eq:NLSDP} 
	and $\bar{\zeta}\in {\cal M}(\bar{x})$ \eqref{eq:defcalm}.
	Then we have the following two results: 
	\begin{itemize}
		\item[(a)] If $\bar{\zeta}\in{\rm ri}\,{\cal M}(\bar{x})$ (the relative interior of ${\cal M}(\bar{x})$), the SOSC \eqref{eq:soscr} holds at $(\bar{x},\bar{\zeta})$ if and only if there are positive constants $\rho_3$, $\theta$, $\varepsilon$, $l$ such that for all $\rho\geq\rho_3$ and all $\zeta\in{\cal M}(\bar{x})\cap\mathbb{B}_{\varepsilon}(\bar{\zeta})$  the uniform quadratic growth condition
		\begin{equation}\label{ieq:ung}
			\mathscr{L}(x,\zeta,\rho)\geq f(\bar{x})+l\|x-\bar{x}\|^2\quad{\rm for\; all}\;x\in\mathbb{B}_{\theta}(\bar{x}) 
		\end{equation}
		is satisfied. 
		\item[(b)] If $\bar{\zeta}\in{\rm rbd}\,{\cal M}(\bar{x})$ (the relative boundary of ${\cal M}(\bar{x})$), the SOSC holds at $(\bar{x},\bar{\zeta})$ if and only if there are positive constants $\rho_3$, $\theta$, $\varepsilon$, $l$ such that \eqref{ieq:ung} holds uniformly for all $\rho\geq\rho_3$ and all $\zeta\in{\cal M}(\bar{x})\cap\mathbb{B}_{\varepsilon}(\bar{\zeta})$ with $\pi(\Gamma)=\pi(\overline{\Gamma})$. 
	\end{itemize}
\end{theorem}
{\bf Proof.} ``$\Leftarrow$" can be obtained from \cite[Theorem 8.4]{mohammadi20}. Then we are going to verify the opposite direction. It follows from Lemma \ref{lemma:soscegeq0} that there exist the positive constants $l'$, $\varepsilon'$ and $\rho_3$ for which \eqref{eq:sosceqine} is satisfied for all $\zeta\in{\cal M}(\bar{x})\cap\mathbb{B}_{\varepsilon'}(\bar{\zeta})$ and all $\rho\geq\rho_3$. Using this and $f(\bar{x})=\mathscr{L}(\bar{x},\zeta,\rho_3)$ for any $\zeta\in{\cal M}(\bar{x})$, which can be obtained by the same proof of \cite[Proposition 3.2 (a)]{HSarabi20}, we deduce from \eqref{eq:sosceqine} that for any given $\zeta\in{\cal M}(\bar{x})\cap\mathbb{B}_{\varepsilon'}(\bar{\zeta})$ there exists $\theta_{\zeta}>0$ for which we have 
\begin{equation}\label{eq:dep}
	\mathscr{L}(x,\zeta,\rho_3)\geq f(\bar{x})+\frac{l'}{2}\|x-\bar{x}\|^2\quad{\rm for\;all}\;x\in\mathbb{B}_{\theta_{\zeta}}(\bar{x}),
\end{equation}
where the constant $l'$ can be chosen the same for all the multipliers $\zeta\in{\cal M}(\bar{x})\cap\mathbb{B}_{\varepsilon'}(\bar{\zeta})$. It can be obtained directly from the definition of the second subderivative. The radii of the balls centered at $\bar{x}$ in \eqref{eq:dep}, however, depend on $\zeta$. If $\bar{\zeta}\in{\rm ri}\,{\cal M}(\bar{x})$, we argue that a common radius can be chosen for all the multipliers $\zeta\in{\cal M}(\bar{x})$ that are sufficiently close to $\bar{\zeta}$. 
Its proof is exactly the same as the proof of \cite[Theorem 4.5]{HMSarabi}. We omit it here for simplicity.  

If $\bar{\zeta}\in{\rm rbd}\,{\cal M}(\bar{x})$, we prove  that a common radius can be chosen for all $\zeta\in{\cal M}(\bar{x})\cap\mathbb{B}_{\varepsilon}(\bar{\zeta})$ with $\pi(\Gamma)=\pi(\overline{\Gamma})$. Following the proof of \cite[Proposition 3.2]{HSarabi20}, we know that $f(\bar{x})=\mathscr{L}(\bar{x},\zeta,\rho_3)$. It is easy to see that for all $x\in\mathbb{B}_{\theta_{\bar{\zeta}}}(\bar{x})$, we have 
\begin{equation*}
	\frac{f(\bar{x})-\mathscr{L}(x,\bar{\zeta},\rho_3)}{\|x-\bar{x}\|^2}\leq-\frac{l'}{2}.
\end{equation*}
From Proposition \ref{prop:uni-soel} with $\zeta\in{\cal M}(\bar{x})$ and the positive homogenous of degree 2 of second semiderivative, we have 
$$\frac{f(\bar{x})-\mathscr{L}(x,\zeta,\rho_3)}{\|x-\bar{x}\|^2}=-\frac{1}{2}{\rm d}_x^2\mathscr{L}(\bar{x},\zeta,\rho_3)(\frac{x-\bar{x}}{\|x-\bar{x}\|})+{O}(\|x-\bar{x}\|),$$
where ${O}(\|x-\bar{x}\|)$ is uniform for all $\zeta\in{\cal M}(\bar{x})\cap\mathbb{B}_{r}(\bar{\zeta})$ with $\pi(\Gamma)=\pi(\overline{\Gamma})$, with its uniform radius $c$ and uniform constant $q_1$. 

From the proof of the Lemma \ref{lemma:soscegeq0}, we know for all $\zeta\in{\cal M}(\bar{x})\cap\mathbb{B}_{\varepsilon'}(\bar{\zeta})$, 
$${\rm d}_x^2\mathscr{L}(\bar{x},\zeta,\rho_3)(\frac{x-\bar{x}}{\|x-\bar{x}\|})\geq{\rm d}_x^2\mathscr{L}(\bar{x},\bar{\zeta},\rho_3)(\frac{x-\bar{x}}{\|x-\bar{x}\|})+O(\|\Gamma-\overline{\Gamma}\|).$$
Suppose the uniform constant for $O(\|\Gamma-\overline{\Gamma}\|)$ is $q_2$. Let $\|x-\bar{x}\|\leq \min\{\theta_{\bar{\lambda}},c,l'/(16q_1)\}:=\theta$ and $\|\zeta-\bar{\zeta}\|\leq \min\{r,\varepsilon',l'/(16q_2)\}:=\varepsilon$ with $\zeta\in{\cal M}(\bar{x})$, $\pi(\Gamma)=\pi(\overline{\Gamma})$. It is easy to see that $\|\Gamma-\overline{\Gamma}\|\leq\|\zeta-\bar{\zeta}\|\leq\varepsilon$. Thus we have 
\begin{align*}
	&\frac{f(\bar{x})-\mathscr{L}(x,\zeta,\rho_3)}{\|x-\bar{x}\|^2}=-\frac{1}{2}{\rm d}_x^2\mathscr{L}(\bar{x},\zeta,\rho_3)(\frac{x-\bar{x}}{\|x-\bar{x}\|})+{O}(\|x-\bar{x}\|)\\
	&\leq-\frac{1}{2}{\rm d}_x^2\mathscr{L}(\bar{x},\bar{\zeta},\rho_3)(\frac{x-\bar{x}}{\|x-\bar{x}\|})+O(\|\Gamma-\overline{\Gamma}\|)+\frac{l'}{16}\\
	&\leq-\frac{1}{2}{\rm d}_x^2\mathscr{L}(\bar{x},\bar{\zeta},\rho_3)(\frac{x-\bar{x}}{\|x-\bar{x}\|})+\frac{l'}{8}\\
	&=\frac{f(\bar{x})-\mathscr{L}(x,\bar{\zeta},\rho_3)}{\|x-\bar{x}\|^2}+{O}(\|x-\bar{x}\|)+\frac{l'}{8}\leq\frac{f(\bar{x})-\mathscr{L}(x,\bar{\zeta},\rho_3)}{\|x-\bar{x}\|^2}+\frac{3l'}{16}.
\end{align*}
Taking the supremum of $x\in\mathbb{B}_{\theta}(\bar{x})$ on both sides and we have for all $\zeta\in\mathbb{B}_{\varepsilon}(\bar{\zeta})\cap{\cal M}(\bar{x})$ with $\pi(\Gamma)=\pi(\overline{\Gamma})$, 
$$\sup\limits_{x\in\mathbb{B}_{\theta}(\bar{x})}\frac{f(\bar{x})-\mathscr{L}(x,\lambda,\rho_3)}{\|x-\bar{x}\|^2}\leq\sup\limits_{x\in\mathbb{B}_{\theta}(\bar{x})}\frac{f(\bar{x})-\mathscr{L}(x,\bar{\lambda},\rho_3)}{\|x-\bar{x}\|^2}+\frac{3l'}{16}\leq-\frac{5l'}{16}.$$
It follows that for all $x\in\mathbb{B}_{\theta}(\bar{x})$ and $\zeta\in\mathbb{B}_{\varepsilon}(\bar{\zeta})\cap{\cal M}(\bar{x})$ with $\pi(\Gamma)=\pi(\overline{\Gamma})$, 
$$\mathscr{L}(x,\zeta,\rho_3)\geq f(\bar{x})+\frac{5l'}{16}\|x-\bar{x}\|^2. $$
From \cite[Exercise 11.56]{RoWe98}, we know that for all $\rho\geq\rho_3$, $\mathscr{L}(x,\zeta,\rho)\geq\mathscr{L}(x,\zeta,\rho_3)$. Setting $l=\frac{5l'}{16}$ and we have proved \eqref{ieq:ung}. 
$\BOX$

\vskip 10pt
The locally Lipschitz continuous property of local minimizer of augmented Lagrangian function for NLSDP is established in \cite[Theorem 1]{SSZhang08} under nondegeneracy and strong SOSC. It is worth noting that \cite[Proposition 5.2]{HSarabi20} have verified the uniformly isolated calmness of the local minimizers of the augmented Lagrangian for composite linear quadratic problems by only requiring SOSC. Here, we aim at extending a similar result to NLSDP, which may relax the conditions in \cite[Theorem 1]{SSZhang08}. 
\begin{proposition}\label{lemma:uniform-conv}
	Let $\bar{x}\in{\cal X}$ be a stationary point to the NLSDP \eqref{eq:NLSDP} 
	and $\bar{\zeta}\in {\cal M}(\bar{x})$ \eqref{eq:defcalm}.
	Suppose $(\bar{x},\bar{\zeta})$ satisfies SOSC \eqref{eq:soscr}. Then there are positive constants $\tau, \rho_3, \widehat{r}>0$ such that for every $\rho\geq\rho_3$ and $\zeta\in\mathbb{B}_{\widehat{r}/2\tau}(\bar{\zeta})$, 
	the set of the local minimizers of function ${\cal L}(x,\zeta,\rho)$ over $x\in\mathbb{B}_{\widehat{r}}(\bar{x})$, defined by ${\cal S}_{\rho}(\zeta)$, satisfies the uniform isolated calmness property, i.e., 
	$${\cal S}_{\rho}(\zeta)\subseteq\{\bar{x}\}+\tau\|\zeta-\bar{\zeta}\|\mathbb{B}$$
	and satisfies $\varnothing\neq {\cal S}_{\rho}(\zeta)\subseteq{\rm int}\,\mathbb{B}_{\widehat{r}}(\bar{x})$, where $\mathbb{B}$ is the unite ball in ${\cal X}$. 
\end{proposition}
{\bf Proof.} 
It follows from the continuity of $\mathscr{L}$ and the compactness of $\mathbb{B}_{\widehat{r}}(\bar{x})$ that ${\cal S}_{\rho}(\zeta)\neq\varnothing$ (cf. \cite[Theorem 4.16]{Rudin}). From \cite{RoWe98}, we know that the Lagrangian function $\mathscr{L}(x,\zeta,\rho)$ is concave on $\zeta$. Combined with \cite[Theorem 8.4]{mohammadi20}, we have there are $\rho_3>0$, $\widehat{r}>0$ and $l>0$ such that for all $\rho\geq\rho_3$ and $x\in\mathbb{B}_{\widehat{r}}(\bar{x})$ 
\begin{align*}
	\mathscr{L}(x,\zeta,\rho)&\geq\mathscr{L}(x,\bar{\zeta},\rho)-\langle\nabla_{\zeta}\mathscr{L}(x,\zeta,\rho),\bar{\zeta}-\zeta\rangle\\
	&=\mathscr{L}(x,\bar{\zeta},\rho)-\langle\Phi(x)-\Pi_{\cal K}(\Phi(x)+\rho^{-1}\zeta),\bar{\zeta}-\zeta\rangle\\
	&\geq f(\bar{x})+l\|x-\bar{x}\|^2-\langle\Phi(x)-\Pi_{\cal K}(\Phi(x)+\rho^{-1}\zeta),\bar{\zeta}-\zeta\rangle, 
\end{align*}
where the $l,\rho_3$ and $\widehat{r}$ is the same as in \cite[Theorem 8.4]{mohammadi20},  ${\cal K}=\{0\}\times{\cal S}^n_+$. 
Let $u\in{\cal S}_{\rho}(\zeta):=\arg\min\{\mathscr{L}(x,\zeta,\rho)\mid x\in\mathbb{B}_{\widehat{r}}(\bar{x})\}$ and we have 
$$\mathscr{L}(u,\zeta,\rho)\leq\mathscr{L}(\bar{x},\zeta,\rho)=f(\bar{x})+\frac{\rho}{2}{\rm dist}(\Phi(\bar{x})+\rho^{-1}\zeta,{\cal K})^2-\frac{1}{2}\rho^{-1}\|\zeta\|^2\leq f(\bar{x}).$$
Define $\tau:=\frac{\kappa_{\Phi}}{l}+\sqrt{\frac{\kappa_{\Phi}^{2}}{l^{2}}+\frac{1}{l\rho_3}}$ and fix $\zeta\in\mathbb{B}_{\widehat{r}/2\tau}(\bar{\zeta})$ and $\rho\geq\rho_3$, where $\kappa_{\Phi}$ is the Lipschitzian constant of $\Phi$. 
Combining the above two inequalities together, we have 
\begin{equation}\label{eq:impun}
	\|u-\bar{x}\|^2\leq\frac{1}{l}\langle\Phi(u)-\Pi_{\cal K}(\Phi(u)+\rho^{-1}\zeta),\bar{\zeta}-\zeta\rangle.
\end{equation}
Since $\bar{\zeta}\in{N}_{{\cal K}}(\Phi(\bar{x}))$, we have $\Phi(\bar{x})=\Pi_{{\cal K}}(\Phi(\bar{x})+\rho^{-1}\bar{\zeta})$. It follows that 
\begin{align*}
	&\|\Phi(u)-\Pi_{\cal K}(\Phi(u)+\rho^{-1}\zeta)\|\\
	&=\|\Phi(u)-\Phi(\bar{x})+\Pi_{{\cal K}}(\Phi(\bar{x})+\rho^{-1}\bar{\zeta})-\Pi_{\cal K}(\Phi(u)+\rho^{-1}\zeta)\|\\
	&\leq2\|\Phi(u)-\Phi(\bar{x})\|+\rho^{-1}\|\bar{\zeta}-\zeta\|\leq2\kappa_{\Phi}\|u-\bar{x}\|+\rho^{-1}\|\bar{\zeta}-\zeta\|.  
\end{align*}
Combining the above two inequalities together leads us to 
$$\|u-\bar{x}\|^2\leq\frac{1}{l}\big(2\kappa_{\Phi}\|u-\bar{x}\|+\rho^{-1}\|\bar{\zeta}-\zeta\|\big)\|\zeta-\bar{\zeta}\|.$$
The latter inequality can be written as the following form 
$$l\|u-\bar{x}\|^2-2\kappa_{\Phi}\|\zeta-\bar{\zeta}\|\cdot\|u-\bar{x}\|-\rho^{-1}\|\zeta-\bar{\zeta}\|^2\leq0. $$
It follows that 
$$\|u-\bar{x}\|\leq\big(\frac{\kappa_{\Phi}}{l}+\sqrt{\frac{\kappa_{\Phi}^{2}}{l^{2}}+\frac{1}{l {\rho}}}\big)\|\zeta-\bar{\zeta}\|\leq\tau\|\zeta-\bar{\zeta}\|\leq\tau\widehat{r}/2\tau<\widehat{r}. $$
Then we have completed the proof. 
$\BOX$
\begin{remark}\label{rem:ime}
	Recalling the definition of residual function $R(x,\zeta)$ \eqref{eq:residual func}. It is easy to know that for KKT point $(\bar{x},\bar{\zeta})$, there exist $r_2$ and $\kappa_2$ such that for all $(x,\zeta)\in\mathbb{B}_{r_2}(\bar{x},\bar{\zeta})$, 
	\begin{equation}\label{eq:invbound}
		R(x,\zeta)\leq\kappa_2\big(\|x-\bar{x}\|+{\rm dist}(\zeta,{\cal M}(\bar{x}))\big). 
	\end{equation}
	Its proof is in the same way as in \cite[Proposition 5.4]{HSarabi20}. 
		Moreover, by Theorem \ref{thm:equuqgc} and the proof of \eqref{eq:impun}, we can prove in the same approach that if $\bar{\zeta}\in{\rm ri}\,{\cal M}(\bar{x})$, for all $\mu\in{\cal M}(\bar{x})\cap\mathbb{B}_{\varepsilon}(\bar{\zeta})$, there exist $\rho_3>0$, $\theta>0$ and $l>0$ such that for all $\rho\geq\rho_3$, $x\in\mathbb{B}_{\theta}(\bar{x})\cap{\cal S}_{\rho}(\zeta)$, $\zeta\in\mathcal{Y}\times{\cal S}^n$,  
	\begin{equation}\label{eq:sosc with qgc}
		\|x-\bar{x}\|^2\leq\frac{1}{l}\langle\Phi(x)-\Pi_{\cal K}(\Phi(x)+\rho^{-1}\zeta),\mu-\zeta\rangle.
	\end{equation}
	Similarly, if $\bar{\zeta}\in{\rm rbd}\,{\cal M}(\bar{x})$, we have that there also exists $\varepsilon>0$ such that for all $\mu=(y,\Gamma)\in{\cal M}(\bar{x})\cap\mathbb{B}_{\varepsilon}(\bar{\zeta})$ with $\pi({\Gamma})=\pi({\overline{\Gamma}})$, \eqref{eq:sosc with qgc} holds for all $\rho\geq\rho_3$, $x\in\mathbb{B}_{\theta}(\bar{x})\cap{\cal S}_{\rho}(\zeta)$, $\zeta\in\mathcal{Y}\times{\cal S}^n$. 
\end{remark}

\begin{remark}
	It is worth noting that \cite[Theorems 1 and  2]{Rockafellarc} obtained augmented tilt stability under the variational sufficient condition. The relationship between augmented tilt stability and uniform isolated calmness of ${\cal S}_{\rho}(\zeta)$ remains unknown to us though it seems the former is stronger. However, as mentioned in \cite{Rockafellarcal}, the variational sufficient condition used in \cite{Rockafellarc} may fail when SOSC holds.  The variational sufficient condition can be satisfied for fully amenable problems, which include NLP, nonlinear second-order cone programming (NLSOC) and exclude NLSDP obviously. This implies that the approach taken in this paper is different from that of \cite{Rockafellarc}. 
\end{remark}

\subsection{Local convergence analysis}
Now we are going to establish the linear convergence of ALM for NLSDP. 
The following error bound estimate is an analogy to \cite[Theorem 5.5]{{HSarabi20}}, which mainly focuses on the polyhedron case. We illustrate it in two different cases. 
\begin{proposition}\label{prop:kk bound}
	Let $\bar{x}\in{\cal X}$ be a stationary point to the NLSDP \eqref{eq:NLSDP} 
	and $\bar{\zeta}\in {\cal M}(\bar{x})$ \eqref{eq:defcalm}.  Suppose $(\bar{x},\bar{\zeta})$ satisfies SOSC  \eqref{eq:soscr} and the semi-isolated calmness (see Definition \ref{def:semi-ic}) holds for $S_{KKT}$ at $\big((0,0),(\bar{x},\bar{\zeta})\big)$. 
	If $\bar{\zeta}\in{\rm ri}\,{\cal M}(\bar{x})$, 
		then there exists positive constants $r_3$, $\kappa_3$ and $\rho_3$ such that for all $\rho\geq\rho_3$, $(x,\zeta)\in\mathbb{B}_{r_3}(\bar{x},\bar{\zeta})$ with $R(x,\zeta)>0$, and all the optimal solutions $u$ to problem
		\begin{equation}\label{eq:local subp}
			\min\mathscr{L}(w,\zeta,\rho)\quad{\rm subject\;to}\quad w\in\mathbb{B}_{\widehat{r}}(\bar{x})
		\end{equation}
		with $\widehat{r}$ obtained in Proposition \ref{lemma:uniform-conv},   
		the error bound estimate
		\begin{equation}\label{eq:alm errorb}
			\|u-x\|+\|\nabla e_{1/\rho}\delta_{\mathcal{K}}(\Phi(u)+\rho^{-1}\zeta)-\zeta\|\leq\kappa_3 R(x,\zeta)
		\end{equation}
		holds, where ${\cal K}=\{0\}\times{\cal S}^n_+$.
		If $\bar{\zeta}\in{\rm rbd}\,{\cal M}(\bar{x})$, \eqref{eq:alm errorb} also holds for all $\rho\geq\rho_3$, $(x,\zeta)\in\mathbb{B}_{r_3}(\bar{x},\bar{\zeta})$ with $\pi(\Gamma_{\pi})=\pi({\overline{\Gamma}})$ and $R(x,\zeta)>0$ , where $\Pi_{{\cal M}(\bar{x})}(\zeta)=(y_{\pi},\Gamma_{\pi})$. 
	\end{proposition}
	{\bf Proof.} By Proposition \ref{lemma:uniform-conv}, we know that for every $\zeta\in\mathbb{B}_{\widehat{r}/2r}(\bar{\zeta})$ and every $\rho\geq\rho_3$ any optimal solution $u$ to \eqref{eq:local subp} satisfies the first-order optimality condition
	\begin{equation}\label{eq:prooffo}
		\nabla_x\mathscr{L}(u,\zeta,\rho)=0.
	\end{equation}
If $\bar{\zeta}\in{\rm ri}\,{\cal M}(\bar{x})$, assume by contradiction that the error bound estimate \eqref{eq:alm errorb} fails, which implies that there exists a sequence $\{(x^k,\zeta^k,\rho^k)\}_{k=1}^{\infty}\subset\mathcal{X}\times\mathcal{Y}\times{\cal S}^n\times[\rho_3,\infty)$ with $(x^k,\zeta^k)\rightarrow(\bar{x},\bar{\zeta})$ as $k\rightarrow\infty$ such that for each $k$,
	\begin{equation}\label{eq:contradiction}
		\|u^k-x^k\|+\|d^k-\zeta^k\|>kR_k, 
	\end{equation}
where $d^k:=\nabla e_{1/\rho_k}\delta_{{\cal K}}(\Phi(u^k)+\frac{\zeta^k}{\rho^k})$ and $u^k$ is an optimal solution to \eqref{eq:local subp} for $(\zeta,\rho)=(\zeta^k,\rho^k)$ 
	and 
	\begin{equation}\label{eq:def-Rk}
		R_k:=R(x^k,\zeta^k)=\|\nabla_xL(x^k,\zeta^k)\|+\|\Phi(x^k)-\Pi_{\cal K}(\Phi(x^k)+\zeta^k)\|.
	\end{equation} 
If $\bar{\zeta}\in{\rm rbd}\,{\cal M}(\bar{x})$, we also assume by contradiction and the only difference lies in the supposed sequence $\zeta^k$ satisfies $\pi(\Gamma_{\pi}^k)=\pi(\overline{\Gamma})$ in addition. 
	For each $k$, denote $c_k:=\|u^k-x^k\|+\|d^k-\zeta^k\|$. Thus, we know from \eqref{eq:contradiction} that $R_k=o(c_k)$. It then follows from \eqref{eq:def-Rk} that for $k\to \infty$
	\begin{equation}\label{eq:proof3}
		\nabla_xL(x^k,\zeta^k)+o(c_k)=0\quad{\rm and}\quad \Phi(x^k)+o(c_k)=\Pi_{{\cal K}}(\Phi(x^k)+\zeta^k).
	\end{equation}
	By passing to a subsequence if necessary, we are able to find $0\neq(\xi,\eta)$ with $\xi\in\mathcal{X}$ and $ \eta:=(\eta_0,\eta_1)\in {\cal Y}\times{\cal S}^n$ such that 
	\begin{equation}\label{eq:proof1}
		\frac{u^k-x^k}{c_k}\rightarrow\xi\quad{\rm and}\quad\frac{d^k-\zeta^k}{c_k}\rightarrow\eta \quad{\rm as}\quad k\to \infty.
	\end{equation}
	Since $(x^k,\zeta^k)\rightarrow(\bar{x},\bar{\zeta})$, we obtain from \eqref{eq:invbound} that $R_k\rightarrow0$ as $k\to \infty$. Moreover, it follows from the definition of semi-isolated calmness and \cite[Theorem 3.1]{KSteck17} that there exists $\kappa>0$ such that for $(x,\zeta)$ sufficiently close to $(\bar{x},\bar{\zeta})$, 
	\begin{equation*}
		\|x-\bar{x}\|+{\rm dist}(\zeta,{\cal M}(\bar{x}))\le \kappa R(x,\zeta).
	\end{equation*}
	For each $k$, let $\mu^k=(\mu^k_1,\mu^k_2)\in{\cal Y}\times{\cal S}^n$ be the metric projection of $\zeta^k$ over the nonempty closed convex set ${\cal M}(\bar{x})$, i.e., $\mu^k:=\Pi_{{\cal M}(\bar{x})}(\zeta^k)$. Thus, without loss of generality, we may assume that $x^k-\bar{x}=O(R_k)$ and $\zeta^k-\mu^k=O(R_k)$ for $k\to \infty$, which in turn results in 
	\begin{equation}\label{eq:proof2}
		x^k-\bar{x}=o(c_k)\quad{\rm and}\quad\zeta^k-\mu^k=o(c_k)\quad{\rm as}\quad k\rightarrow\infty.
	\end{equation} 
	The latter along with $\zeta^k\rightarrow\bar{\zeta}$ tells us that $\mu^k\rightarrow\bar{\zeta}$ as $k\to \infty$, since ${\cal M}(\bar{x})$ is closed and convex. Therefore, we get $\mu^k\in{\cal M}(\bar{x})\cap\mathbb{B}_{\varepsilon}(\bar{\zeta})$ for all $k$ sufficiently large.
	
	If $\bar{\zeta}\in{\rm ri}\,{\cal M}(\bar{x})$, then we know from \eqref{eq:sosc with qgc} in Remark \ref{rem:ime} that there exists $l>0$  for all $k$ sufficiently large, 
	\begin{equation}\label{prop12proof2}
		\|u^k-\bar{x}\|^2\leq\frac{1}{\rho^kl}\langle d^k-\zeta^k,\mu^k-\zeta^k\rangle\leq\frac{1}{\rho^kl}\| d^k-\zeta^k\|\cdot\|\mu^k-\zeta^k\|.
	\end{equation}
	If $\bar{\zeta}\in{\rm rbd}\,{\cal M}(\bar{x})$, again from Remark \ref{rem:ime}, we know that \eqref{prop12proof2} also holds for all $k$ sufficiently large when $\pi(\mu^k_2)=\pi(\overline{\Gamma})$. For simplicity, we only show the $\bar{\zeta}\in{\rm ri}\,{\cal M}(\bar{x})$ case here while the other case can be obtained similarly. 
	It follows from \cite[Theorem 2.26]{RoWe98} that for each $k$,
	\begin{equation}\label{eq:dk_exp}
		d^k=\nabla e_{1/\rho_k}\delta_{{\cal K}}(\Phi(u^k)+\frac{\zeta^k}{\rho^k})=\rho^k\big(\Phi(u^k)+\frac{\zeta^k}{\rho^k}-\Pi_{\cal K}(\Phi(u^k)+\frac{\zeta^k}{\rho^k})\big), 
	\end{equation}
which implies that
	\begin{align}
		&\|d^k-\zeta^k\|=\|\rho^k\big(\Phi(u^k)+\frac{\zeta^k}{\rho^k}-\Pi_{\cal K}(\Phi(u^k)+\frac{\zeta^k}{\rho^k})\big)-\zeta^k\|=\rho^k\|\Phi(u^k)-\Pi_{\cal K}(\Phi(u^k)+\frac{\zeta^k}{\rho^k})\|\nonumber\\
		&\leq\rho^k\big({\rm dist}(\Phi(u^k),{\cal K})+\|\Pi_{\cal K}(\Phi(u^k))-\Pi_{\cal K}(\Phi(u^k)+\frac{\zeta^k}{\rho^k})\|\big)\nonumber\\
		&\leq\rho^k\big({\rm dist}(\Phi(u^k),{\cal K})+\|\frac{\zeta^k}{\rho^k}\|\big)\leq\rho^k\big(\|\Phi(\bar{x})-\Phi(u^k)\|+{\rm dist}(\Phi(\bar{x}),{\cal K})+\|\frac{\zeta^k}{\rho^k}\|\big).\label{prop12proof1}
	\end{align}
	Combining \eqref{prop12proof2}, \eqref{prop12proof1} and the boundedness of $\{u^k\}_{k=1}^{\infty}$ from Proposition \ref{lemma:uniform-conv} together, since $\zeta^k\rightarrow\bar{\zeta}$ and $\mu^k\rightarrow\bar{\zeta}$, we obtain that 
	\begin{align}\label{eq:ukrabx}
		\|u^k-\bar{x}\|^2&\leq\frac{1}{l}\big(\|\Phi(\bar{x})-\Phi(u^k)\|+\|\frac{\zeta^k}{\rho^k}\|\big)\cdot\|\mu^k-\zeta^k\|\rightarrow0\quad{\rm as}\quad k\rightarrow\infty.
	\end{align}
	Multipling \eqref{prop12proof2} by $1/c_k^2$, we know from \eqref{eq:proof1} and \eqref{eq:proof2} that 
	\begin{equation}\label{eq:proof4}
		\frac{\|u^k-\bar{x}\|^2}{c_k^2}\leq\frac{1}{\rho^kl}\frac{\| d^k-\zeta^k\|}{c_k}\cdot\frac{\| \mu^k-\zeta^k\|}{c_k}\rightarrow0\quad{\rm as}\quad k\rightarrow\infty,
	\end{equation} 
	which implies for $k\to \infty$, $u^k-\bar{x}=o(c_k)$. This, together with \eqref{eq:proof2}, yields that 
	\begin{equation}\label{eq:xir0}
		\xi=\lim\limits_{k\rightarrow\infty}\frac{u^k-x^k}{c_k}=\lim\limits_{k\rightarrow\infty}\frac{u^k-\bar{x}}{c_k}-\lim\limits_{k\rightarrow\infty}\frac{x^k-\bar{x}}{c_k}=0-0=0.
	\end{equation}
	
	The desired result then follows by contradiction, if we show that $\eta=0$. To this end, let us consider the following two cases.
	
{\bf Case 1:} either the sequence $\{\rho^k\}_{k=1}^{\infty}$ or $\{\rho^k/c_k\}_{k=1}^{\infty}$ is bounded, it follows from \eqref{eq:proof4} or \eqref{eq:ukrabx} that $\frac{\rho^k}{c_k}\|u^k-\bar{x}\|\rightarrow0$ as $k\rightarrow\infty$. It follows from \eqref{eq:dk_exp} and the twice continuous differentiability of $\Phi$ that there exists $\kappa_{\Phi}>0$ such that 
	\begin{align*}
		\frac{\|d^k-\zeta^k\|}{c_k}&=\frac{\rho^k}{c_k}\|\Phi(u^k)-\Pi_{\cal K}(\Phi(u^k)+(\rho^k)^{-1}\zeta^k)\|\\
		&=\frac{\rho^k}{c_k}\|\Phi(u^k)-\Phi(\bar{x})+\Pi_{{\cal K}}(\Phi(\bar{x})+(\rho^k)^{-1}\mu^k)-\Pi_{\cal K}(\Phi(u^k)+(\rho^k)^{-1}\zeta^k)\|\\
		&\leq\frac{\rho^k}{c_k}(2\|\Phi(u^k)-\Phi(\bar{x})\|+(\rho^k)^{-1}\|\mu^k-\zeta^k\|)\\
		&\leq2\kappa_{\Phi}\frac{\rho^k}{c_k}\|u^k-\overline{x}\|+\frac{\|\mu^k-\zeta^k\|}{c_k}\rightarrow0\quad{\rm as}\quad k\rightarrow\infty.
	\end{align*}
	Therefore, by \eqref{eq:proof2} we obtain that $\eta=0$. 
	
{\bf Case 2:} both sequences $\{\rho^k\}_{k=1}^{\infty}$ and $\{\rho^k/c_k\}_{k=1}^{\infty}$ are unbounded. By passing to a subsequence if necessary, we may assume that 
	\begin{equation}\label{eq:proof5}
		\rho^k\rightarrow\infty\quad{\rm and}\quad\frac{\rho^k}{c_k}\rightarrow\infty\quad{\rm as}\quad k\rightarrow\infty.
	\end{equation}
	Since $u^k$ is an optimal solution to \eqref{eq:local subp} associated with $(\zeta^k,\rho^k)$, we deduce from \eqref{eq:prooffo} that for each $k$,
	\begin{equation}\label{eq:defzkeq}
	\nabla_x\mathscr{L}(u^k,\zeta^k,\rho^k)=\nabla_x L(u^k,d^k)=0.
	\end{equation}
	From  \eqref{eq:proof4}, we have 
	\begin{equation}\label{eq:12leq0}
	\frac{\rho^k}{c_k}\|u^k-\bar{x}\|^2\leq\frac{1}{l}\frac{\| d^k-\zeta^k\|}{c_k}\cdot\|\mu^k-\zeta^k\|\rightarrow0\quad{\rm as}\quad k\rightarrow\infty. 
	\end{equation}

	For each $k$, define $z^k:=\Phi(u^k)+(\rho^k)^{-1}(\zeta^k-d^k)$. Thus, it follows from \eqref{eq:dk_exp} and \cite[Theorem 2.26]{RoWe98} that $d^k\in\partial \delta_{\mathcal{K}}(z^k)$. Using \eqref{eq:proof1}, \eqref{eq:ukrabx} and \eqref{eq:proof5} allows us to arrive at
	\begin{equation*}\label{eq:defzk}
	z^k=\Phi(u^k)-\frac{d^k-\zeta^k}{c_k}\cdot\frac{c_k}{\rho^k}\rightarrow\Phi(\bar{x})=:\bar{z}\quad{\rm as}\quad k\rightarrow\infty.
	\end{equation*}
	 For each $k$, by noting that $\mu^k\in{\cal M}(\bar{x})$, we have $\mu^k\in{N}_{{\cal K}}(\bar{z})$. 
By using the firmly non-expansiveness of projection mapping (cf. \cite[(1.6)]{z1971}) and \cite[Theorem 31.5]{Rockafellar1970}, we have for each $k$,
$$\|z^k-\bar{z}\|^2\leq\langle z^k-\bar{z},z^k+d^k-\bar{z}-\mu^k\rangle=\|z^k-\bar{z}\|^2+\langle z^k-\bar{z},d^k-\mu^k\rangle,$$
which yields $\langle z^k-\bar{z},d^k-\mu^k\rangle\geq0$. Therefore, we have for each $k$,
$$\frac{\rho^k}{c_k^2}\langle z^k-\bar{z},\mu^k-d^k\rangle=\langle\frac{\rho^k}{c_k}(\Phi(u^k)-\Phi(\bar{x}))-\frac{d^k-\zeta^k}{c_k},\frac{\mu^k-d^k}{c_k}\rangle\leq0, $$
which implies
\begin{equation}\label{eq:prto0}
\langle-\frac{d^k-\zeta^k}{c_k},\frac{\mu^k-d^k}{c_k}\rangle\leq-\langle\frac{\rho^k}{c_k}(\Phi(u^k)-\Phi(\bar{x})),\frac{\mu^k-d^k}{c_k}\rangle.
\end{equation}
Meanwhile, since $\Phi$ is twice continuously differentiable and \eqref{eq:ukrabx} holds, we have $\Phi(u^k)-\Phi(\bar{x})=\nabla\Phi(\bar{x})(u^k-\bar{x})+O(\|u^k-\bar{x}\|^2)$ and $\Phi(u^k)-\Phi(\bar{x})=\nabla\Phi(u^k)(u^k-\bar{x})+O(\|u^k-\bar{x}\|^2)$ as $k\to \infty$. It is worth to note the second $O(\|u^k-\bar{x}\|^2)$  is uniform for $k$ sufficiently large. 
Thus, since $f$ is twice continuous differentiable, we know that there exists $\kappa_f>0$ such that for $k$ sufficiently large,
\begin{align}
&-\langle\frac{\rho^k}{c_k}(\Phi(u^k)-\Phi(\bar{x})),\frac{\mu^k-d^k}{c_k}\rangle=-\frac{\rho^k}{c_k^2}\big(\langle\Phi(u^k)-\Phi(\bar{x}), \mu^k\rangle-\langle\Phi(u^k)-\Phi(\bar{x}), d^k\rangle\big)\nonumber\\
&=-\frac{\rho^k}{c_k^2}\big(\langle \nabla\Phi(\bar{x})(u^k-\bar{x})+O(\|u^k-\bar{x}\|^2), \mu^k\rangle-\langle\nabla\Phi(u^k)(u^k-\bar{x})+O(\|u^k-\bar{x}\|^2), d^k\rangle\big)\nonumber\\
&=-\frac{\rho^k}{c_k^2}\big(\langle (u^k-\bar{x}), \nabla\Phi(\bar{x})^*\mu^k\rangle-\langle(u^k-\bar{x}), \nabla\Phi(u^k)^*d^k\rangle+\langle O(\|u^k-\bar{x}\|^2), \mu^k\rangle+\langle O(\|u^k-\bar{x}\|^2), d^k\rangle\big)\nonumber\\
&=-\frac{\rho^k}{c_k^2}\big(\langle (u^k-\bar{x}), \nabla\Phi(\bar{x})^*\mu^k-\nabla\Phi(u^k)^*d^k\rangle+\langle O(\|u^k-\bar{x}\|^2), \mu^k\rangle+\langle O(\|u^k-\bar{x}\|^2), d^k\rangle\big)\nonumber\\
&=-\frac{\rho^k}{c_k^2}\big(\langle u^k-\bar{x}, \nabla f(u^k)-\nabla f(\bar{x})\rangle+\langle O(\|u^k-\bar{x}\|^2), \mu^k\rangle+\langle O(\|u^k-\bar{x}\|^2), d^k\rangle\big)\nonumber\\
&\leq\frac{\rho^k}{c_k^2}\kappa_f\|u^k-\bar{x}\|^2+\langle O(\frac{\rho^k}{c_k^2}\|u^k-\bar{x}\|^2), -\mu^k\rangle+\langle O(\frac{\rho^k}{c_k}\|u^k-\bar{x}\|^2), -\frac{d^k}{c_k}\rangle\nonumber\\
&=\frac{\rho^k}{c_k^2}\kappa_f\|u^k-\bar{x}\|^2+\langle O(\frac{\rho^k}{c_k^2}\|u^k-\bar{x}\|^2), -\mu^k\rangle+\langle O(\frac{\rho^k}{c_k}\|u^k-\bar{x}\|^2), \frac{-d^k+\zeta^k}{c_k}\rangle+\langle O(\frac{\rho^k}{c_k^2}\|u^k-\bar{x}\|^2), -\zeta^k\rangle, \label{eq:crtcorr}
\end{align}
where the forth equation follows from $\nabla_x L(u^k,d^k)=0$ \eqref{eq:defzkeq} and $\nabla_x L(\bar{x},\mu^k)=0$ for all $k$. It then follows from \eqref{eq:proof4} and \eqref{eq:12leq0} that $$\frac{\rho^k}{c_k^2}\|u^k-\bar{x}\|^2\rightarrow0 \quad \mbox{and} \quad \frac{\rho^k}{c_k}\|u^k-\bar{x}\|^2\rightarrow0 \quad \mbox{as}\quad k\to \infty.$$ 
Combining this with $\zeta^k\rightarrow\bar{\zeta}$, $\mu^k\rightarrow\bar{\zeta}$, \eqref{eq:proof1} and 
\begin{equation*}
\lim_{k\rightarrow\infty}\frac{\mu^k-d^k}{c_k}=\lim_{k\rightarrow\infty}\frac{\mu^k-\zeta^k+\zeta^k-d^k}{c_k}=\lim_{k\rightarrow\infty}\frac{\mu^k-\zeta^k}{c_k}+\lim_{k\rightarrow\infty}\frac{\zeta^k-d^k}{c_k}=-\eta,
\end{equation*}
we obtain from  \eqref{eq:prto0} and \eqref{eq:crtcorr} that $\|\eta\|^2\leq0$ by taking $k\rightarrow\infty$. Thus, we know that $\eta=0$, which completes the proof. 
	$\BOX$

	\vskip 10pt
	Before we put forward the main result of this paper, we need to propose the following assumption, which is needed in the proof of our main result. 
	\begin{assumption}\label{ass1}
	Given $\bar{\zeta}\in{\rm rbd}\,{\cal M}(\bar{x})$. There exits $r_5>0$, $\chi>0$ such that for all $(x,\zeta)\in\mathbb{B}_{r_5}(\bar{x},\bar{\zeta})$ with $\zeta=(y,\Gamma)\notin{\cal M}(\bar{x})$, there exist $\widehat{\zeta}\in{\cal M}(\bar{x})$ with $\pi(\widehat{\Gamma})=\pi(\overline{\Gamma})$ such that 
		\begin{equation*}
			\|\Pi_{{\cal M}(\bar{x})}(\zeta)-\widehat{\zeta}\|\leq \chi R(x,\zeta). 
		\end{equation*}
	\end{assumption}
It is worth to note that Assumption \ref{ass1} is not like any constraint qualification that we are familiar with. However, it is easy to see that it holds at least in the following circumstances. 
\begin{itemize}
\item[(a)] If for all $\Gamma_1$, $\Gamma_2\in{\cal M}(\bar{x})$, $\pi(\Gamma_1)=\pi(\Gamma_2)$ holds, then Assumption \ref{ass1} trivially holds under semi-isolated calmness for $S_{KKT}$ at $\big((0,0),(\bar{x},\bar{\zeta})\big)$ as we only need to pick $\widehat{\zeta}=\Pi_{{\cal M}(\bar{x})}(\zeta)$.  This circumstance also includes the case where ${\cal M}(\bar{x})$ is a singleton. It is worth to note that for non-polyheral cases, the strict Ronbinson condition mentioned in \cite[Definition 2.3]{KSteck19} is not equivalent to the uniqueness of ${\cal M}(\bar{x})$ (cf. \cite[Example 3]{DSZhang17}). 
\item[(b)] Suppose $A=G(\bar{x})+\overline{\Gamma}$ possesses the eigenvalue decomposition \eqref{eq:eig-decomp}. If there is no multiple root in $\beta\cup\gamma$ part eigenvalues of $\overline{\Gamma}$, then for all $\lambda$ sufficiently close to $\bar{\lambda}$, every $\Gamma$ also do not possess multiple root in $\beta\cup\gamma$ part for eigenvalues. We can also pick $\widehat{\zeta}=\Pi_{{\cal M}(\bar{x})}(\zeta)$ under the semi-isolated calmness for $S_{KKT}$ at $\big((0,0),(\bar{x},\bar{\zeta})\big)$ to show the validity of Assumption \ref{ass1}. A sufficient condition for semi-isolated calmness will be given in Section \ref{sec:scsic}. It is worth to note that in this case we only require the existence of such $\overline{\Gamma}$ and an oracle to find a starting point lies in the neighborhood mentioned in Assumption \ref{ass1} when \eqref{eq:soscr} holds for all $\zeta\in{\cal M}(\bar{x})$), which is called Robinson's SOSC \cite{Robinson82} (see also \cite[(25)]{DSZhang17}).  Moreover, in this case, neither the set ${\cal M}(\bar{x})$ has to be a polyhedron nor the strict complementarity condition, i.e., $\overline{\zeta}\in{\rm ri}\,{\cal M}(\bar{x})$ is satisfied.  
\end{itemize}
Moreover, In Section \ref{sec:scsic}, we give 2 examples (Example \ref{example:ass1hold}, \ref{example:ass1hold2}) to show the validity of this assumption. 
	By adding Assumption \ref{ass1} to the conditions in Proposition \ref{prop:kk bound}, we can get a further result of Proposition \ref{prop:kk bound} for the case when $\bar{\zeta}\in{\rm rbd}\,{\cal M}(\bar{x})$. 
	\begin{corollary}\label{coro-u}
		Besides the conditions in Proposition \ref{prop:kk bound}, if the Assumption \ref{ass1} also holds, 
		we have when $\bar{\zeta}\in{\rm rbd}\,{\cal M}(\bar{x})$, \eqref{eq:alm errorb} also holds for all $\rho\geq\rho_3$ and $(x,\zeta)\in\mathbb{B}_{r_3}(\bar{x},\bar{\zeta})$ with $R(x,\zeta)>0$ and $R(\bar{x},\zeta)>0$. 
	\end{corollary}
	{\bf Proof.} The proof is exactly the same as the proof of Proposition \ref{prop:kk bound}. Suppose the contradiction sequence $(x^k,\zeta^k)$ also satisfies $\zeta^k\notin{\cal M}(\bar{x})$. 
	The only difference lies in \eqref{eq:proof2} and \eqref{prop12proof2}, as we can find $\widehat{\mu}^k\in{\cal M}(\bar{x})$ with $\pi(\widehat{\Gamma}^k)=\pi(\overline{\Gamma})$ satisfies $\|\widehat{\mu}^k-\mu^k\|=O(R_k)$. Then $\mu^k$ in \eqref{eq:proof2} and \eqref{prop12proof2} can be alternated by $\widehat{\mu}^k$. 
	Then we have completed the proof. 
	$\BOX$

	\begin{remark}\label{rem:inexact}
		Given $(x^k,\zeta^k)\in\mathbb{B}_{r_3}(\bar{x},\bar{\zeta})$ and $\rho^k\geq\rho_3$ with $r_3$ and $\rho_3$ taken from Proposition \ref{prop:kk bound}. Suppose $x^{k+1}$ is the optimal solution to \eqref{eq:local subp}. Increasing $\kappa_3$ if necessary. 
		Similarly as in \cite[Remark 5.6]{HSarabi20}, we know that for $\widetilde{x}^{k+1}$ sufficiently close to $x^{k+1}$, 
		we also have 
		$$\|\nabla_x\mathscr{L}(\widetilde{x}^{k+1},\zeta^k,\rho^k)\|\leq\epsilon_k$$
		and
		$$\|\widetilde{x}^{k+1}-x^k\|+\|\nabla e_{1/\rho}\delta_{\mathcal{K}}(\Phi(\widetilde{x}^{k+1})+(\rho^k)^{-1}\zeta^k)-\zeta^k\|\leq\kappa_3R(x^k,\zeta^k).$$
	\end{remark}
	
	Next we are going to propose the main result of this paper, which is inspired by \cite[Theorem 5.6]{HSarabi20}. It illustrates the local linear convergence of ALM for NLSDP without requiring the uniqueness of multipliers by applying Proposition \ref{prop:kk bound} and \eqref{eq:sosc with qgc}. 
	\begin{theorem}\label{thm:mainconv}
		Let $\bar{x}\in{\cal X}$ be a stationary point to the NLSDP \eqref{eq:NLSDP} 
		and $\bar{\zeta}\in {\cal M}(\bar{x})$ \eqref{eq:defcalm}. Suppose $(\bar{x},\bar{\zeta})$ satisfies SOSC  \eqref{eq:soscr} and the semi-isolated calmness (see Definition \ref{def:semi-ic}) holds for $S_{KKT}$ at $\big((0,0),(\bar{x},\bar{\zeta})\big)$. 
		\begin{itemize}
			\item[(i)] If $\bar{\zeta}\in{\rm ri}\,{\cal M}(\bar{x})$, then there exist positive constants $\overline{r}$, $\overline{\sigma}$, $\overline{\varrho}$ such that for any starting point $(x^0,\zeta^0)\in\mathbb{B}_{\overline{r}}(\bar{x},\bar{\zeta})$  the primal-dual sequence $\{(x^k,\zeta^k)\}_{k\geq0}$ generated by Algorithm \ref{algo1} with $\rho^k\geq\overline{\varrho}$ and $\epsilon_k=o(R(x^k,\zeta^k))$ for all $k$ satisfies the estimate
			\begin{equation}\label{eq:estn}
				\|x^{k+1}-x^k\|+\|\zeta^{k+1}-\zeta^k\|\leq\overline{\sigma}R(x^k,\zeta^k).
			\end{equation}
		\item[(ii)] If $\bar{\zeta}\in{\rm rbd}\,{\cal M}(\bar{x})$ and Assumption \ref{ass1} holds, then there exist positive constants $\overline{r}$, $\overline{\sigma}$, $\overline{\varrho}$ such that for any starting point $(x^0,\zeta^0)\in\mathbb{B}_{\overline{r}}(\bar{x},\bar{\zeta})$  the primal-dual sequence $\{(x^k,\zeta^k)\}_{k\geq0}$ generated by Algorithm \ref{algo1} with $\rho^k\geq\overline{\varrho}$ and $\epsilon_k=o(R(x^k,\zeta^k))$ and $\zeta^k\notin {\cal M}(\bar{x})$ for all $k$ satisfies the estimate \eqref{eq:estn}.
		\end{itemize}		
	Moreover, for each case, the sequence is convergent to $(\bar{x},\widetilde{\zeta})$ for some $\widetilde{\zeta}\in{\cal M}(\bar{x})$ and its rate of convergence is linear, i.e., for $k$ sufficiently large, 
	\begin{equation}\label{eq:convinb}
		\|(x^{k+1},\zeta^{k+1})-(\bar{x},\widetilde{\zeta})\|\leq \tau^k\|(x^{k},\zeta^{k})-(\bar{x},\widetilde{\zeta})\|,
	\end{equation}
	where $\tau^k=2\sqrt{2}\overline{\sigma}\kappa_1\kappa_2^2(R_k^{-1}\epsilon_k+(\rho^k)^{-1}\overline{\sigma})$, $R_k:=R(x^k,\zeta^k)$. 
\end{theorem}
{\bf Proof.} Consider $R_k:=R(x^k,\zeta^k)$. If $R_k=0$ for some $k$, then the pair $(x^k,\zeta^k)$ satisfies the KKT system and the algorithm should stop. Thus we assume $R_k>0$ for all $k\in\mathbb{N}$. Pick $\kappa_1$ and $r_1$ from Definition \ref{def:semi-ic} with $\mathbb{V}=\mathbb{B}_{r_1}(\bar{x},\bar{\zeta})$ and $\kappa=\kappa_1$, $\kappa_2$ and $r_2$ from \eqref{eq:invbound}, $\rho_3$, $\kappa_3$ and $r_3$ from Proposition \ref{prop:kk bound} or Corollary \ref{coro-u}, $\tau$ and $\widehat{r}$ from Proposition \ref{lemma:uniform-conv}. By the definition of $\epsilon_k$, we can find $r_4>0$ such that 
\begin{equation}\label{thm5eq7}
	\epsilon(x,\zeta)\leq\frac{1}{4\kappa_1\kappa_2}R(x,\zeta)\quad{\rm whenever}\quad(x,\zeta)\in\mathbb{B}_{r_4}(\bar{x},\bar{\zeta}).
\end{equation}
Define $\overline{\sigma}=\kappa_3$ and 
\begin{equation}\label{eqr}
	\overline{r}=\frac{r'}{1+2\sqrt{2}\,\overline{\sigma}\kappa_2}\quad{\rm with}\quad r'=\min\{\widehat{r},\frac{\widehat{r}}{2\tau}, \frac{r_1}{\sqrt{2}\,\overline{\sigma}\kappa_2+1}, r_2, r_4,r_3\}. 
\end{equation}
Pick $q\in(0,1)$ and $\overline{\rho}=\max\{\rho_3,\frac{2\sqrt{2}\kappa_1\kappa_2^2\overline{\sigma}^2}{q}, 4\kappa_1\kappa_2\overline{\sigma}\}$. 
By induction, we want to show that if $\bar{\zeta}\in{\rm ri}\,{\cal M}(\bar{x})$, for any starting point $(x^0,\zeta^0)\in\mathbb{B}_{\overline{r}}(\bar{x},\bar{\zeta})$ the sequence generated by the algorithm with $\rho^k\geq\overline{\rho}$ and $\epsilon_k=o(R(x^k,\zeta^k))$, we have for all $k=0,1,\dots$ the following relationships 
\begin{align}
	&(x^k,\zeta^k)\in\mathbb{B}_{r'}(\bar{x},\bar{\zeta}),\label{ver1}\\
	&\|\nabla_xL(x^{k+1},\zeta^{k+1})\|\leq\epsilon_k,\label{ver2}\\
	&\|x^{k+1}-x^k\|+\|\zeta^{k+1}-\zeta^k\|\leq\overline{\sigma}R_k\label{ver3}
\end{align}
hold. By induction, we firstly assume $k=0$. Since $(x^0,\zeta^0)\in\mathbb{B}_{\overline{r}}(\bar{x},\bar{\zeta})$ and $\overline{r}\leq r'$, we have \eqref{ver1} holds. By Proposition \ref{lemma:uniform-conv} and $\zeta^0\in\mathbb{B}_{\widehat{r}/2\tau}(\bar{\zeta})$, we can find $\widehat{x}^1\in{\rm int}\,\mathbb{B}_{\widehat{r}}(\bar{x})$ satisfying $\nabla_x\mathscr{L}(\widehat{x}^1,\zeta^0,\rho^0)=0$. From Remark \ref{rem:inexact}, we know that we can find $x^1$ sufficiently close to $\widehat{x}^1$ such that $x^1$ satisfies the two relationships in Remark \ref{rem:inexact}.  Define further $\zeta^1=\rho^0[\Phi(x^1)+\zeta^0/\rho^0-\Pi_{\cal K}(\Phi(x^1)+\zeta^0/\rho^0)]$ and we have 
$$\|\nabla_xL(x^1,\zeta^1)\|=\|\nabla_x\mathscr{L}(x^1,\zeta^0,\rho^0)\|\leq\epsilon_k.$$
It follows from Proposition \ref{prop:kk bound} and Remark \ref{rem:inexact} that 
$$\|x^1-x^0\|+\|\zeta^1-\zeta^0\|\leq\overline{\sigma}R_0.$$
Thus, $(x^1,\zeta^1)$ is well defined and satisfies \eqref{ver2} and \eqref{ver3} for $k=0$. Then we assume $(x^k,\zeta^k)$, $k=0,1,\dots,s+1$ are well defined and \eqref{ver1}-\eqref{ver3} hold for $k=0,1,\dots,s$. We now verify the existence of $(x^{s+2},\zeta^{s+2})$ and \eqref{ver1}-\eqref{ver3} satisfies for $k=s+1$. We first show that $(x^{s+1},\zeta^{s+1})\in\mathbb{B}_{r'}(\bar{x},\bar{\zeta})$. Fix an integer $k$ with $0\leq k\leq s$. Since $(x^k,\zeta^k)\in\mathbb{B}_{r'}(\bar{x},\bar{\zeta})$, it follows from \eqref{eq:invbound} that
\begin{align}
	R_k&\leq\kappa_2(\|x^k-\bar{x}\|+{\rm dist}(\zeta^k,{\cal M}(\bar{x})))\label{proof5:eq1}\\
	&\leq\kappa_2(\|x^k-\bar{x}\|+\|\zeta^k-\bar{\zeta}\|)\leq\sqrt{2}\kappa_2\|(x^k,\zeta^k)-(\bar{x},\bar{\zeta})\|\leq\sqrt{2}\kappa r'.\nonumber
\end{align}
Thus we have 
\begin{align*}
	\|(x^{k+1},\zeta^{k+1})-(\bar{x},\bar{\zeta})\|&\leq\|x^{k+1}-x^k\|+\|\zeta^{k+1}-\zeta^k\|+\|(x^k,\zeta^k)-(\bar{x},\bar{\zeta})\|\\
	&\leq\overline{\sigma}R_k+r'\leq(\sqrt{2}\,\overline{\sigma}\kappa_2+1)r'\leq r_1,
\end{align*}
which implies that $(x^{k+1},\zeta^{k+1})\in\mathbb{B}_{r_1}(\bar{x},\bar{\zeta})$. From Definition \ref{def:semi-ic} combined with \cite[Theorem 3.1]{KSteck17}, we obtain
\begin{align*}
	&\|x^{k+1}-\bar{x}\|+{\rm dist}(\zeta^{k+1},{\cal M}(\bar{x}))\\
	&\leq\kappa_1R_{k+1}=\kappa_1(\|\nabla_xL(x^{k+1},\zeta^{k+1})\|+\|\Phi(x^{k+1})-\Pi_{\cal K}(\Phi(x^{k+1})+\zeta^{k+1})\|)\\
	&\leq\kappa_1\epsilon_{k}+\kappa_1\|\Phi(x^{k+1})-\Pi_{\cal K}(\Phi(x^{k+1})+\zeta^{k+1})\|.
\end{align*}
Let $p^{k+1}=\Pi_{\cal K}(\Phi(x^{k+1})+(\rho^k)^{-1}\zeta^k)$. From Algorithm \ref{algo1} we have $\Phi(x^{k+1})-p^{k+1}=(\rho^k)^{-1}(\zeta^{k+1}-\zeta^k)$. It follows from $\zeta^{k+1}=\nabla e_{1/\rho_k}\delta_{\mathcal{K}}(\Phi(x^{k+1})+(\rho^k)^{-1}\zeta^k)$ and $\nabla e_rg(x)=(rI+(\partial g)^{-1})^{-1}(x)$ (see \cite[Theorem 2.26]{RoWe98}) that $\zeta^{k+1}\in{ N}_{\cal K}(p^{k+1})$. By the nonexpansivness of $y\mapsto y-\Pi_{\cal K}(y+\zeta^{k+1})$, we obtained
\begin{align*}
	&\|\Phi(x^{k+1})-\Pi_{\cal K}(\Phi(x^{k+1})+\zeta^{k+1})\|\nonumber\\
	&=\|\Phi(x^{k+1})-\Pi_{\cal K}(\Phi(x^{k+1})+\zeta^{k+1})\|-\|p^{k+1}-\Pi_{\cal K}(p^{k+1}+\zeta^{k+1})\|\nonumber\\
	&\leq\|\Phi(x^{k+1})-\Pi_{\cal K}(\Phi(x^{k+1})+\zeta^{k+1})-(p^{k+1}-\Pi_{\cal K}(p^{k+1}+\zeta^{k+1}))\|\nonumber\\
	&\leq\|\Phi(x^{k+1})-p^{k+1}\|. 
\end{align*}
Thus we have 
$$\|x^{k+1}-\bar{x}\|+{\rm dist}(\zeta^{k},{\cal M}(\bar{x}))\leq\kappa_1\epsilon_{k}+(\rho^k)^{-1}\kappa_1\|\zeta^{k+1}-\zeta^k\|,$$
which can be further calculated as 
\begin{align} 
	\|x^{k+1}-\bar{x}\|+{\rm dist}(\zeta^{k+1},{\cal M}(\bar{x}))&\overset{\eqref{ver3}}{\leq}\kappa_1\epsilon_k+\frac{\overline{\sigma}\kappa_1}{\rho_k}R_k\overset{\eqref{thm5eq7}}{\leq}\frac{1}{4\kappa_2}R_k+\frac{1}{4\kappa_2}R_k\label{th5eq6}\\
	&\overset{\eqref{eq:invbound}}{\leq}\frac{1}{2}\big(\|x^{k}-\bar{x}\|+{\rm dist}(\zeta^{k},{\cal M}(\bar{x}))\big).\nonumber
\end{align}
It follows that 
\begin{equation}\label{proof2th5}
	\|x^{k+1}-\bar{x}\|+{\rm dist}(\zeta^{k+1},{\cal M}(\bar{x}))\leq\frac{1}{2^{k+1}}\big(\|x^{0}-\bar{x}\|+{\rm dist}(\zeta^{0},{\cal M}(\bar{x}))\big).
\end{equation}
Then we have 
\begin{align*}
	&\|(x^{s+1},\zeta^{s+1})-(x^0,\zeta^0)\|\leq\sum\limits_{k=0}^s\|(x^{k+1},\zeta^{k+1})-(x^k,\zeta^k)\|\overset{\eqref{ver3}}{\leq}\overline{\sigma}\sum\limits_{k=0}^sR_k\\
	&\overset{\eqref{proof5:eq1}}{\leq}\overline{\sigma}\kappa_2\sum\limits_{k=0}^s\big(\|x^{k}-\bar{x}\|+{\rm dist}(\zeta^{k},{\cal M}(\bar{x}))\big)\\
	&\overset{\eqref{proof2th5}}{\leq}\overline{\sigma}\kappa_2\sum\limits_{k=0}^s\frac{1}{2^k}\big(\|x^{0}-\bar{x}\|+{\rm dist}(\zeta^{0},{\cal M}(\bar{x}))\big)\\
	&\leq2\overline{\sigma}\kappa_2\big(\|x^{0}-\bar{x}\|+{\rm dist}(\zeta^{0},{\cal M}(\bar{x}))\big)\leq2\overline{\sigma}\kappa_2\big(\|x^{0}-\bar{x}\|+\|\zeta^0-\bar{\zeta}\|\big).
\end{align*}
Thus we arrive at the estimate
\begin{align*}
	&\|(x^{s+1},\zeta^{s+1})-(\bar{x},\bar{\zeta})\|\leq\|(x^{s+1},\zeta^{s+1})-(x^0,\zeta^0)\|+\|(x^0,\zeta^0)-(\bar{x},\bar{\zeta})\|\\
	&\leq2\overline{\sigma}\kappa_2(\|x^0-\bar{x}\|+\|\zeta^0-\bar{\zeta}\|)+\|(x^0,\zeta^0)-(\bar{x},\bar{\zeta})\|\\
	&\leq(2\sqrt{2}\overline{\sigma}\kappa_2+1)\|(x^0,\zeta^0)-(\bar{x},\bar{\zeta})\|\leq(2\sqrt{2}\overline{\sigma}\kappa_2+1)\overline{r}=r',
\end{align*}
where the last inequality comes from $(x^0,\zeta^0)\in\mathbb{B}_{\overline{r}}(\bar{x},\bar{\zeta})$. Then we have verified $(x^{s+1},\zeta^{s+1})\in\mathbb{B}_{r'}(\bar{x},\bar{\zeta})$. By \eqref{eqr}, we get $\zeta^{s+1}\in\mathbb{B}_{\widehat{r}/2\tau}(\bar{\zeta})$, and hence Proposition \ref{lemma:uniform-conv} ensures the optimal solution $\widehat{x}^{s+2}$ such that $\widehat{x}^{s+2}\in{\rm int}\,\mathbb{B}_{\widehat{r}}(\bar{x})$. Thus we have $\nabla_x\mathscr{L}(\widehat{x}^{s+2},\zeta^{s+1},\rho^{s+1})=0$. Still from Remark \ref{rem:inexact}, we can find $x^{s+2}$ sufficiently close to $\widehat{x}^{s+2}$ such that $x^{s+2}$ satisfies the two relationships in Remark \ref{rem:inexact}  
and we observe that 
$$\|\nabla_xL(x^{s+2},\zeta^{s+2})\|=\|\nabla_x\mathscr{L}(x^{s+2},\zeta^{s+1},\rho^{s+1})\|\leq\epsilon_{s+1}.$$
By Proposition \ref{prop:kk bound} and Remark \ref{rem:inexact}, we have 
$$\|x^{s+2}-x^{s+1}\|+\|\zeta^{s+2}-\zeta^{s+1}\|\leq\overline{\sigma}R_{s+1}.$$
Then we have finished verifying \eqref{ver1}-\eqref{ver3} for $k=s+1$. If $\bar{\zeta}\in{\rm rbd}\,{\cal M}(\bar{x})$ and Assumption \ref{ass1} holds, we want to prove for any starting point $(x^0,\zeta^0)\in\mathbb{B}_{\overline{r}}(\bar{x},\bar{\zeta})$ the sequence generated by the algorithm with $\rho^k\geq\overline{\rho}$, $\epsilon_k=o(R(x^k,\zeta^k))$ and $R(\bar{x},\zeta^k)>0$ for all $k$, relationships \eqref{ver1}-\eqref{ver3} also hold for all $k=0,1,\dots$. The proof is exactly the same as the $\bar{\zeta}\in{\rm ri}\,{\cal M}(\bar{x})$ case and the only difference lies in the Proposition \ref{prop:kk bound} used above should be alternated by Corollary \ref{coro-u}.

Then we prove the convergence of the sequence. Use the same argument as in the proofs of \eqref{proof2th5}, we have 
\begin{equation}\label{eq5th5}
	\|(x^{k+l},\zeta^{k+l})-(x^k,\zeta^k)\|\leq2\overline{\sigma}\kappa_2\big(\|x^k-\overline{x}\|+{\rm dist}(\zeta^k,{\cal M}(\bar{x}))\big)\quad{\rm for\; all}\quad k,l\in\mathbb{N}.
\end{equation}
It follows from \eqref{proof2th5} that the righthand side of \eqref{eq5th5} goes to $0$ as $k\rightarrow\infty$, which implies $\{(x^k,\zeta^k)\}$ is Cauchy. Assume $\{(x^k,\zeta^k)\}$ convergences to $(\bar{x},\widetilde{\zeta})$, where $\widetilde{\zeta}\in{\cal M}(\bar{x})$. Let $l\rightarrow\infty$ in \eqref{eq5th5}. Then we have 
$$\|(x^k,\zeta^k)-(\bar{x},\widetilde{\zeta})\|\leq2\overline{\sigma}\kappa_2\big(\|x^k-\bar{x}\|+{\rm dist}(\zeta^k,{\cal M}(\bar{x}))\big),$$
which together with \eqref{th5eq6} verifies
\begin{align*}
	&\|(x^{k+1},\zeta^{k+1})-(\bar{x},\widetilde{\zeta})\|\leq2\overline{\sigma}\kappa_2\big(\|x^{k+1}-\bar{x}\|+{\rm dist}(\zeta^{k+1},{\cal M}(\bar{x}))\big)\\
	&\leq2\overline{\sigma}\kappa_2\kappa_1(R_k^{-1}\epsilon_k+(\rho^k)^{-1}\overline{\sigma})R_k\\
	&\leq2\overline{\sigma}\kappa_2^2\kappa_1(R_k^{-1}\epsilon_k+(\rho^k)^{-1}\overline{\sigma})\big(\|x^k-\bar{x}\|+{\rm dist}(\zeta^k,{\cal M}(\bar{x}))\big)\\
	&\leq2\sqrt{2}\overline{\sigma}\kappa_1\kappa_2^2(R_k^{-1}\epsilon_k+(\rho^k)^{-1}\overline{\sigma})\|(x^k,\zeta^k)-(\bar{x},\widehat{\zeta})\|.
\end{align*}
Combining this with $\rho^k\geq\overline{\rho}$ and $\epsilon_k=o(R_k)$ result in 
$$\lim\sup\limits_{k\rightarrow\infty}\frac{\|(x^{k+1},\zeta^{k+1})-(\bar{x},\widetilde{\zeta})\|}{\|(x^{k},\zeta^{k})-(\bar{x},\widetilde{\zeta})\|}\leq\lim\sup\limits_{k\rightarrow\infty}2\sqrt{2}\overline{\sigma}\kappa_1\kappa_2^2(R_k^{-1}\epsilon_k+(\rho^k)^{-1}\overline{\sigma})\leq q. $$
It follows from $q\in(0,1)$ that the convergence rate is linear. 
Then we have completed the whole proof. 
$\BOX$

\begin{remark}
	If $\rho^k\rightarrow\infty$, we obtain the asymptotic Q-superlinear convergence rate of KKT pair from \eqref{eq:convinb} as $\tau^k\rightarrow0$. Regarding to the update of $\rho^k$, in \cite{KSteck17,HSarabi20}, they apply a practical rule in Algorithm \ref{algo1}, step 4  to update $\rho^k$, i.e., defining the auxiliary function by
	$$V(x,\zeta,\rho):=\|\nabla_x\mathscr{L}(x,\zeta,\rho)\|+\|\Phi(x)-\Pi_{\cal K}(\Phi(x)+\rho^{-1}\zeta)\|.$$ 
Suppose $c\in(0,1)$.	If $k=0$ or $V(x^{k+1},\zeta^k,\rho^k)\leq cV(x^k,\zeta^{k-1},\rho^{k-1})$ holds, set $\rho^{k+1}:=\rho^k$; otherwise, set $\rho^{k+1}:=\varsigma\rho^k$. By taking advantage of Theorem \ref{thm:mainconv} and similar manners as in \cite{HSarabi20}, we are able to obtain the boundedness of  $\{\rho^k\}_{k=1}^{\infty}$. 
	
	In addiction, it is worth to note that when ${\cal M}(\bar{x})$ is a singleton, the semi-isolated calmness of $S_{KKT}$ is reduced to its isolated calmness. As mention in \cite[Theorem 24]{DSZhang17}, the isolated calmness of $S_{KKT}$ is equivalent to SOSC and SRCQ \cite[(4.125)]{bonnans} under RCQ \cite[Definition 2.86]{bonnans}. Thus, the result obtained in Theorem \ref{thm:mainconv} reduces to \cite[Theorem 4.2]{KSteck19} for NLSDP. 
\end{remark}
\begin{remark} In Theorem \ref{thm:mainconv}, if $\bar{\zeta}\in{\rm rbd}\,{\cal M}(\bar{x})$ and Assumption \ref{ass1} holds, we have shown that the primal-dual sequence $\{(x^k,\zeta^k)\}_{k\geq0}$ generated by Algorithm \ref{algo1} converge to $(\bar{x},\widetilde{\zeta})$ for some $\widetilde{\zeta}\in{\cal M}(\bar{x})$  if $\zeta^k\notin{\cal M}(\bar{x})$ for sufficiently large $k$. It remains unknown from our approach whether the primal sequence $\{x^k\}_{k\geq0}$ converges to $\bar{x}$ linearly if the dual sequence $\{\zeta^k\}_{k\geq0}$ terminates finitely. This is a future work we are working on. 
\end{remark}

\section{A sufficient condition for the semi-isolated calmness of $S_{KKT}$ \eqref{eq:KKT-NLSDP-p1}}\label{sec:scsic}
In this section, we give a sufficient condition for the semi-isolated calmness of KKT pair. 
In order to reach the goal, we need the definition of bounded linear regularity of a collection of closed convex sets, which can be found in, e.g., \cite[Definition 5.6]{BBorwein96}. 
\begin{definition}
	Let $D_{1}, D_{2}, \ldots, D_{m} \subseteq \mathcal{X}$ be closed convex sets for some positive integer $m .$ Suppose that $D:=D_{1} \cap D_{2} \cap \ldots \cap D_{m}$ is non-empty. The collection $\left\{D_{1}, D_{2}, \ldots, D_{m}\right\}$ is said to be boundedly linearly regular if for every bounded set $\mathcal{B} \subseteq \mathcal{X}$, there exists a constant $\kappa>0$ such that
	$$
	\operatorname{dist}(x, D) \leqslant \kappa \max \left\{\operatorname{dist}\left(x, D_{1}\right), \ldots, \operatorname{dist}\left(x, D_{m}\right)\right\}, \forall x \in \mathcal{B}.
	$$
\end{definition}
A sufficient condition to guarantee the property of bounded linear regularity is established in \cite[Corollary 3]{BBLi99}. 
Denote ${\cal G}_1(\bar{x})=\{(y,\Gamma)\in\mathcal{Y}\times{\cal S}^n\mid\nabla f(\bar{x})+\nabla h(\bar{x})^*y+\nabla G(\bar{x})^*\Gamma=0\}$ and ${\cal G}_2(\bar{x})=\{(y,\Gamma)\in\mathcal{Y}\times{\cal S}^n\mid\Gamma\in{\cal N}_{{\cal S}_+^n}(G(\bar{x}))\}$. It is easy to see that ${\cal G}_1(\bar{x})$ is a polyhedron and ${\cal G}_2(\bar{x})$ is convex. Along with \cite[Theorem 3.1]{KSteck17}, the following result gives a sufficient condition for semi-isolated calmness. Its proof is inspired from \cite[Theorem 5.9]{MSarabi19} by using a similar reduction method and we omit it here for brevity. 
\begin{theorem}\label{thm: error bound}
	Let $\bar{x}\in{\cal X}$ be a stationary point to the NLSDP \eqref{eq:NLSDPp1} with $({a}_1,{a}_2,{b})=(0,0,0)$ and $(\bar{y},\overline{\Gamma})\in {\cal M}(\bar{x})$. Suppose SOSC  \eqref{eq:soscr} holds at $(\bar{x},\bar{y},\overline{\Gamma})$ and 
	\begin{equation}\label{boundedeq}
		{\cal G}_1(\bar{x})\cap{\rm ri}\,{\cal G}_2(\bar{x})\neq\varnothing. 
	\end{equation}
	Then there exist a constant $\kappa_1>0$, neighborhoods $\mathbb{V}:=\mathbb{B}_{r_1}(\bar{x},\bar{y},\overline{\Gamma})$ of $(\bar{x},\bar{y},\overline{\Gamma})$ and $\mathbb{U}$ of $(0,0,0)$ such that for any $(a_1,a_2,b)\in\mathbb{U}$,
	\begin{equation*}
		\|x-\bar{x}\|+{\rm dist}((y,\Gamma),{\cal M}(\bar{x}))\le \kappa_1\|(a_1,a_2,b)\|\quad \forall\, (x,y,\Gamma)\in{S}_{\rm KKT}(a_1,a_2,b)\cap \mathbb{V}.
	\end{equation*}
\end{theorem}
 Regarding to condition \eqref{boundedeq}, it suffice to admit a point $(y',\Gamma')\in{\cal M}(\bar{x})$ possessing the strict complementarity condition with respect to SDP cone. We have to emphasis that this $(y',\Gamma')$ can be different from both the reference point $(\bar{y},\overline{\Gamma})$ and the limit point of the sequence generated by ALM, which is mentioned in Theorem \ref{thm:mainconv}.

We can also study the validity of semi-isolated calmness directly from \cite[Theorem 5.9]{MSarabi19}. Suppose $(\bar{x},\bar{y},\overline{\Gamma})$ is a KKT pair. It follows that the key of the validity of semi-isolated calmness mainly lies in when ${S}$ is calm at $(\bar{y},\overline{\Gamma})$ for $({a},{b})=(0,0)$, where $S(a,b)=\{(y,\Gamma)\in\mathcal{Y}\times{\cal S}_+^n\mid \nabla_x L(\bar{x},y,\Gamma)=a,\;\Gamma\in{\cal N}_{{\cal S}_+^n}(G(\bar{x})+b)\}$.  
It is easy to see that the strong regularity/isolated calmness/Aubin of $S$ implies the calmness of $S$. (The definition of calmness, isolated calmness and Aubin can be seen in e.g., \cite{YeYe,RoWe98} and \cite{DoRo14}).

To end this paper, we propose two examples to illustrate that the conditions required in Theorem \ref{thm:mainconv} can be satisfied indeed. 
\begin{example}\label{example:ass1hold}
	\begin{equation*}
		\begin{array}{cl}
			\min & \frac{1}{2}x^3\\
			{\rm s.t} & -x^2\left[\begin{array}{ccc}
				0 & 0 & 0  \\ 
				0 & 0 & 0  \\ 
				0 & 0 & 1
			\end{array}\right]\in{\cal S}^3_+, \quad \Leftarrow \Gamma
		\end{array}
	\end{equation*}
	It is easy to see that the optimal solution is $\bar{x}=0$, its corresponding multiplier set is ${\cal M}(\bar{x})=\{\Gamma\mid\Gamma\in{\cal S}_-^3\}$.  
	Pick $\overline{\Gamma}={\rm Diag}(0,-1,-2)$. Then for all $\Gamma\in\mathbb{B}_{{\rm min}\{1/3,r_1\}}(\overline{\Gamma})\backslash{\cal M}(\bar{x})$ with $r_1$ taken from Definition \ref{def:semi-ic} with $\mathbb{V}=\mathbb{B}_{r_1}(\bar{x},\bar{\lambda})$, we know that $$\Pi_{{\cal M}(\bar{x})}(\Gamma)=Q{\rm Diag}(\min\{0,\Gamma_1\},\Gamma_2,\Gamma_3)Q^T,$$ where $Q\in{\cal O}^3(\Gamma)$. Let $\widehat{\Gamma}=Q{\rm Diag}(0,\Gamma_2,\Gamma_3)Q^T$. It follows from Definition \ref{def:semi-ic} and \cite[Theorem 3.1]{KSteck17} that 
	$$\|\Pi_{{\cal M}(\bar{x})}(\Gamma)-\widehat{\Gamma}\|\leq{\rm dist}(\Gamma,{\cal M}(\bar{x}))=O(R(x,\Gamma)).$$
	Thus Assumption \ref{ass1} is satisfied. It is easy to calculate $\nabla^2_{xx}L(\bar{x},\overline{\Gamma})=2>0$, which implies SOSC holds at $(\bar{x},\overline{\Gamma})$. Furthermore, it is obvious that bounded linear regular holds. Then we get the semi-isolated calmness of ${S}_{KKT}$ at $\big(0,(\bar{x},\overline{\Gamma})\big)$ holds by Theorem \ref{thm: error bound}. Thus, if we have an oracle of finding a starting point in $\mathbb{B}_{\overline{r}}(\bar{x},\overline{\Gamma})$ (cf. Theorem \ref{thm:mainconv}), the sequence $\{(x^k,\Gamma^k)\}$ generated by ALM converges to $(\bar{x},\widetilde{\Gamma})$ for some $\widetilde{\Gamma}\in{\cal M}(\bar{x})$ while  $\widetilde{\Gamma}$ does not have to be $\overline{\Gamma}$ or $\widehat{\Gamma}$. 
\end{example}

We also  provide another nontrivial NLSDP example, which is modified from the example proposed in the arxiv version of \cite[Example 2]{CSToh2016} for different purpose.  
\begin{example}\label{example:ass1hold2}
	Consider the following example 
	\begin{equation*}
		\begin{array}{cl}
			\min & \frac{1}{2}x^2+2t\\
			{\rm s.t} & tA-x^2I_2\in{\cal S}^2_+, \quad \Leftarrow \Gamma\\ [3pt]
			&  t\geq0, \quad\quad\quad\quad\quad \Leftarrow y
		\end{array}
	\end{equation*}
	where $A=\left[\begin{array}{cc}
		1 & -2  \\ 
		-2 & 1
	\end{array}\right]$. This problem possesses the unique optimal solution $(\bar{t},\bar{x})=(0,0)$. The corresponding multiplier is 
	$${\cal M}(\bar{t},\bar{x})=\{(\Gamma,y)\in{\cal S}_-^2\times\Re\mid\langle A,-\Gamma\rangle\leq 2\}.$$ 
	We can pick $\bar{y}=0$ and $\overline{\Gamma}=\left[\begin{array}{cc}
		0 & 0  \\ 
		0 & -1
	\end{array}\right]$. For all $\Gamma\in\mathbb{B}_{\min\{r_1, 1/(2\sqrt{10})\}}(\overline{\Gamma})$ with $r_1$ taken from Definition \ref{def:semi-ic} with $\mathbb{V}=\mathbb{B}_{r_1}(\bar{x},\bar{\lambda})$, we know that $\langle A,-\Gamma\rangle<2$, which implies that $\Pi_{{\cal M}(\bar{t},\bar{x})}(\Gamma)=\Pi_{{\cal S}_-^2}(\Gamma)$. Suppose $\Gamma=Q{\rm Diag}(\Gamma_1,\Gamma_2)Q^T$ with $Q\in{\cal O}^2(\Gamma)$. Let $\widehat{\Gamma}=Q{\rm Diag}(0,\Gamma_2)Q^T$. Then we have 
	$$\|\Pi_{{\cal M}(\bar{x})}(\Gamma)-\widehat{\Gamma}\|\leq{\rm dist}(\Gamma,{\cal M}(\bar{x}))=O(R(x,\Gamma)),$$
	which verifies Assumption \ref{ass1}. 
	It is easy to see that $\left[\begin{array}{cc}
		-1/2 & 0  \\ 
		0 & -1/2
	\end{array}\right]\in{\cal G}_1(\bar{t},\bar{x})\cap{\rm ri}\,{\cal G}_2(\bar{t},\bar{x})$, which implies the validity of boundedly linear regularity by \cite[Corollary 3]{BBLi99}. It is easy to check that SOSC holds since $\nabla^2 L(\bar{t},\bar{x},\overline{\Gamma},\bar{y})=\left[\begin{array}{cc}
		0 & 0  \\ 
		0 & 3
	\end{array}\right]$ and for all $w=(w_1,w_2)\in{\cal C}(\bar{t},\bar{x})$, we have $w_1=0$. Thus we have SOSC holds at $(\bar{t},\bar{x},\overline{\Gamma},\bar{y})$. 
	It follows that semi-isolated calmness of ${S}_{KKT}$ at $\big(0,(\bar{t},\bar{x},\overline{\Gamma},\bar{y})\big)$ holds by Theorem \ref{thm: error bound}. 
In this case, if we have an oracle of finding a starting point in $\mathbb{B}_{\overline{r}}(\bar{t},\bar{x},\overline{\Gamma},\bar{y})$ (cf. Theorem \ref{thm:mainconv}), the sequence $\{(t^k,x^k,\Gamma^k,y^k)\}$ generated by ALM converges to $(\bar{t},\bar{x},\widetilde{\Gamma},\widetilde{y})$ for some $(\widetilde{\Gamma},\widetilde{y})\in{\cal M}(\bar{t},\bar{x})$ while  $(\widetilde{\Gamma},\widetilde{y})$ does not have to be $(\overline{\Gamma},\bar{y})$ or $(\widehat{\Gamma},\widehat{y})$. 
\end{example}

\section{Conclusions}\label{sec6}
In this paper, we have shown that the augmented Lagrangian method convergences linearly for NLSDP under certain conditions without requiring the uniqueness of the Lagrangian multiplier. During the establishment of ALM convergence, we obtain the uniform second expansion of the Moreau envelope of SDP and give several sufficient conditions for the semi-isolated calmness of $S_{KKT}$. In the future, we will mainly focus on the following three ongoing works. Firstly, as \cite{CSToh2016} shows, usually, the dual Q-linear convergence rate together with the KKT residual R-linear rate is enough in practical solvers. Therefore it is meaningful to study whether we can get a convergence result of that kind instead of a primal-dual type under a weaker condition. Inspired by work \cite{Rockafellarcal}, which extends the convex framework of ALM convergence to the non-convex case by variational sufficiency, we see hope in extending it to non-convex non-polyhedral problems. 
Secondly, although we consider solving the ALM subproblem inexactly, we have not put forward a practical relative error criterion for it. This is a future work we focus on and \cite{ESilva13,CSToh2016} may provide some inspirations. Thirdly, we are also working on providing sufficient (necessary) conditions to Assumption \ref{ass1}.

\begin{appendices}
{\normalsize 
\section{Proof of Lemma \ref{lemma:unih}}\label{app:A}
For each $s\in\{1,\dots,\bar{d}\}$, let $\Lambda_{\bar{\iota}_s\bar{\iota}_s}={\rm Diag}(\lambda_{\bar{\iota}_s}(A))$ and $\Xi_{\bar{\iota}_s\bar{\iota}_s}={\rm Diag}(\lambda_{\bar{\iota}_s}(A+H))$. We first show that \eqref{eq:unih1} holds. If $\bar{d}=1$, i.e., $\lambda_1(\overline{A})=\dots=\lambda_n(\overline{A})$, the first equation in \eqref{eq:unih1} trivially holds. Next we assume that $\bar{d}\geq2$. From \eqref{eq:unih2}, we have for any ${\cal S}^n\ni H\rightarrow0$, 
\begin{equation*}
\left[\begin{array}{cccc}
\Lambda_{\bar{\iota}_1\bar{\iota}_1} & 0 & \dots & 0\\
0 & \Lambda_{\bar{\iota}_2\bar{\iota}_2} & \dots & 0\\
\vdots&\vdots&\ddots&\vdots\\
0 & 0 &\dots & \Lambda_{\bar{\iota}_{\bar{d}}\bar{\iota}_{\bar{d}}}
\end{array}\right] U+HU=U\left[\begin{array}{cccc}
\Xi_{\bar{\iota}_1\bar{\iota}_1} & 0 & \dots & 0\\
0 & \Xi_{\bar{\iota}_2\bar{\iota}_2} & \dots & 0\\
\vdots&\vdots&\ddots&\vdots\\
0 & 0 &\dots & \Xi_{\bar{\iota}_{\bar{d}}\bar{\iota}_{\bar{d}}}
\end{array}\right].
\end{equation*}
It follows from $A\in\mathbb{B}_{r}(\overline{A})$ that $\lambda_i(A)\neq\lambda_j(A)$ whenever $i\in\bar{\iota}_t$, $j\in\bar{\iota}_s$ with $s\neq t$. 
It is easy to see that for all $i\in\bar{\iota}_s$, $j\in\bar{\iota}_t$ with $s\neq t$, $U_{ij}=\frac{(HU)_{ij}}{\Lambda_{ii}-\Xi_{jj}}$. Then we have for all $\|H\|\leq\omega:=r/6$, 
$$\frac{\|U_{\bar{\iota}_{s} \bar{\iota}_{t}}\|}{\|H\|}\leq\sum_{i\in\bar{\iota}_s, j\in\bar{\iota}_t}\frac{1}{(\Lambda_{ii}-\Xi_{jj})^2}\leq\sum_{i\in\bar{\iota}_k, j\in\bar{\iota}_l}\frac{1}{(\lvert v_i(\overline{A})-v_j(\overline{A})\rvert-2r-2\omega)^2}:=q,$$
where $\omega$ and $q$ are independent of $A$. 
Hence we obtain that 
$$U_{\bar{\iota}_{s} \bar{\iota}_{t}}={O}(\|H\|) \quad \forall\, 1 \leq s \neq t \leq \bar{d},$$
where ${O}(\|H\|)$ is uniform for all $A\in\mathbb{B}_{r}(\overline{A})$. 
By using the fact that $U$ is orthogonal, we obtain directly that the second equation in \eqref{eq:unih1} holds. In order to prove \eqref{eq:unih3}, we consider the SVD of each $U_{\bar{\iota}_{s} \bar{\iota}_{s}}$, $s=1,\dots,\bar{d}$. Fix $s\in\{1,\dots,\bar{d}\}$. Let $W$ and $V$ be in ${\cal O}^{\lvert\bar{\iota}_s\rvert}$ such that $U_{\bar{\iota}_{s} \bar{\iota}_{s}}=W\Sigma V^T$, where $\Sigma$ is a nonnegative diagonal matrix. From \eqref{eq:unih1}, we obtain that for all $A\in\mathbb{B}_{r}(\overline{A})$, 
$$W\Sigma^2W^T=I_{\lvert\bar{\iota}_{s}\rvert}+{O}(\|H\|^2),$$
which is equivalent to 
$$\Sigma^2=W^TW+{O}(\|H\|^2)=I_{\lvert\bar{\iota}_{s}\rvert}+{O}(\|H\|^2).$$
Since $\Sigma$ is a nonnegative diagonal matrix, we may conclude that 
$$\Sigma={\rm Diag}(1+{O}(\|H\|^2),\dots,1+{O}(\|H\|^2)).$$ 
Therefore, from $U_{\bar{\iota}_{s} \bar{\iota}_{s}}=W\Sigma V^T$, we have $U_{\bar{\iota}_{s} \bar{\iota}_{s}}=WV^T+O(\|H\|^2)$. Since $WV^T\in{\cal O}^{\lvert\bar{\iota}_{s}\rvert}$, we know that for all $A\in\mathbb{B}_{r}(\overline{A})$, \eqref{eq:unih3} holds. Next, we shall show \eqref{eq:unih4} holds. 
For each $s\in\{1,\dots,\bar{d}\}$ by comparing the $s$-th diagonal block of both sides of \eqref{eq:unih2}, we obtain that 
\begin{equation}\label{eq:unih6}
U^T_{\bar{\iota}_{s}}(\Lambda(A)+H)U_{\bar{\iota}_{s}}=\Xi_{\bar{\iota}_s\bar{\iota}_s}.
\end{equation}
Fix $s\in\{1,\dots,\bar{d}\}$. From \eqref{eq:unih1} and \eqref{eq:unih6}, we know that 
\begin{align*}
U_{\bar{\iota}_{s}}^T\Lambda(A)U_{\bar{\iota}_{s}}&=\left[\begin{array}{ccc}
{O}(\|H\|)&U_{\bar{\iota}_{s} \bar{\iota}_{s}}&{O}(\|H\|)\end{array}\right] \left[\begin{array}{ccc}
\Lambda_1(A)&0&0\\
0&\Lambda(A)_{\bar{\iota}_{s} \bar{\iota}_{s}}&0\\
0&0&\Lambda_2(A)
\end{array}\right] \left[\begin{array}{c}
{O}(\|H\|)\\
U_{\bar{\iota}_{s} \bar{\iota}_{s}}\\
{O}(\|H\|)
\end{array}\right]\\
&={O}(\|H\|^2)\Lambda_1(A)+U_{\bar{\iota}_{s} \bar{\iota}_{s}}^T\Lambda(A)_{\bar{\iota}_{k} \bar{\iota}_{k}}U_{\bar{\iota}_{k} \bar{\iota}_{k}}+{O}(\|H\|^2)\Lambda_2(A). 
\end{align*}
It follows that 
$$\Xi_{\bar{\iota}_{k} \bar{\iota}_{k}}-\big({O}(\|H\|^2)\Lambda_1(A)+U_{\bar{\iota}_{k} \bar{\iota}_{s}}^T\Lambda(A)_{\bar{\iota}_{s} \bar{\iota}_{s}}U_{\bar{\iota}_{s} \bar{\iota}_{s}}+{O}(\|H\|^2)\Lambda_2(A)\big)=U_{\bar{\iota}_{s} \bar{\iota}_{s}}^TH_{\bar{\iota}_{s} \bar{\iota}_{s}}U_{\bar{\iota}_{s} \bar{\iota}_{s}}+{O}(\|H\|^2).$$
Since $U_{\bar{\iota}_{s} \bar{\iota}_{s}}=Q_s+{O}(\|H\|^2)$ and $\|\Lambda(A)\|\leq\|\Lambda(\overline{A})\|+r$, we obtain that 
$$Q_s^TH_{\bar{\iota}_{s} \bar{\iota}_{s}}Q_k=\Xi_{\bar{\iota}_{s} \bar{\iota}_{s}}-Q_s^T\Lambda(A)_{\bar{\iota}_{s} \bar{\iota}_{s}}Q_s+{O}(\|H\|^2). $$
Hence \eqref{eq:unih4} holds with the uniform ${O}(\|H\|^2)$ for all $A\in\mathbb{B}_{r}(\overline{A})$. 
The proof is completed. 

\section{Proof of Proposition \ref{prop:unisocpi}}\label{app:B} 
Firstly, we show \eqref{eq:piuni} holds for the case that $A=\Lambda(A)$. For any $H\in{\cal S}^n$, denote $Z=A+H$. Let $U\in{\cal O}^n$ (depending on $H$) be such that 
\begin{equation}\label{eq:unisocpi2}
\Lambda(A)+H=U\Lambda(Z)U^T.
\end{equation}
Let $\omega>0$ be any fixed number such that $0\leq\omega\leq\frac{{\lambda}_{\lvert\alpha\rvert}(\overline{A})-r}{2}$ if $\alpha\neq\varnothing$ and be any fixed positive number otherwise. Then, define the following continuous scalar function 
$$
f(t):=\left\{\begin{array}{ll}
t & \text { if } \quad t>\omega \\
2 t-\omega & \text { if } \quad \frac{\omega}{2}<t<\delta \\
0 & \text { if } \quad t<\frac{\omega}{2}.
\end{array}\right.
$$
Therefore, we have 
 $$
\left\{\lambda_{1}(A), \ldots, \lambda_{\lvert\alpha\rvert}(A)\right\} \in(\omega,+\infty) \quad \text { and } \quad\left\{\lambda_{\lvert\alpha\rvert+1}(A), \ldots, \lambda_{n}(A)\right\} \in(-\infty, \frac{\omega}{2}). 
$$
For the scalar function $f$, let $F:{\cal S}^n\rightarrow{\cal S}^n$  be the corresponding L{\"o}wner’s operator, i.e., for any $W\in{\cal S}^n$, 
$$F(W):=\sum_{i=1}^nf(\lambda_i(W))P_iP_i^T,$$
where $P\in{\cal O}^n(W)$. Since $f$ is real analytic on the open set $(-\infty,\frac{\omega}{2})\cup(\omega,+\infty)$, It is well-known that for $H$ sufficiently close to zero,
 \begin{equation}\label{eq:unisocpi3}
 F(A+H)-F(A)-F'(A)H=O(\|H\|^2)
 \end{equation}
 and 
 $$F'(A)H=\left[\begin{array}{ccc}
H_{\alpha \alpha} & H_{\alpha \beta} & \Sigma_{\alpha \gamma} \circ H_{\alpha \gamma} \\
H_{\alpha \beta}^{T} & 0 & 0 \\
\Sigma_{\alpha \gamma}^{T} \circ H_{\alpha \gamma}^{T} & 0 & 0
\end{array}\right], $$
 where $O(\|H\|^2)$ is independent of $A$ for any $A\in\mathbb{B}_r(\overline{A})$ and $\Sigma\in{\cal S}^n$  is given by 
 $$\Sigma_{ij}=\frac{\max\{\lambda_i(A),0\}-\max\{\lambda_j(A),0\}}{\lambda_i(A)-\lambda_j(A)},\quad i,j=1,\dots,n.$$
Let $R(\cdot):=\Pi_{{\cal S}_+^n}(\cdot)-F(\cdot)$. By the definition of $f$, we know that $F(A)=\Pi_{{\cal S}_+^n}(A)$, which implies that $R(A)=0$.  Meanwhile, it is clear that the matrix valued function $R$ is directionally differentiable at $A$, and the directional derivative of $R$ for any given direction $H\in{\cal S}^n$ is given by
$$
R^{\prime}(A ; H)=\Pi_{\mathcal{S}_{+}^{n}}^{\prime}(A ; H)-F^{\prime}(A) H=\left[\begin{array}{ccc}
0 & 0 & 0 \\
0 & \Pi_{\mathcal{S}_{+}^{\lvert\beta\rvert}}( H_{\beta \beta}) & 0 \\
0 & 0 & 0
\end{array}\right]. 
$$
By the Lipschitz continuity of $\lambda(\cdot)$, we know that for $H$ sufficiently close to zero, 
$$\left\{\lambda_{1}(Z), \ldots, \lambda_{\lvert\alpha\rvert}(Z)\right\} \in(\omega,+\infty) , \quad\left\{\lambda_{\lvert\alpha\rvert+1}(Z), \ldots, \lambda_{\lvert\beta\rvert}(Z)\right\} \in(-\infty, \frac{\omega}{2})
$$
and 
$$\{\lambda_{\lvert\beta\rvert+1}(Z), \dots, \lambda_{n}(Z)\} \in(-\infty,0).$$
Therefore, by the definition of $F$ , we know that for $H$ sufficiently close to zero,
$$R(A+H)=\Pi_{{\cal S}_+^n}(A+H)-F(A+H)=U\left[\begin{array}{ccc}
0 & 0 & 0 \\
0 & \Pi_{\mathcal{S}_{+}^{\lvert\beta\rvert}}(\Lambda(Z)_{\beta \beta}) & 0 \\
0 & 0 & 0
\end{array}\right]U^T.$$
Since $U\in{\cal O}^n(Z)$, we know from Lemma \ref{lemma:unih} that for any ${\cal S}^n\ni H\rightarrow0$, there exists an orthogonal matrix $Q\in{\cal O}^{\lvert\beta\rvert}$ such that 
\begin{equation}\label{eq:unisocpi1}
U_{\beta}=\left[\begin{array}{c}
{O}(\|H\|)\\
U_{\beta\beta}\\
{O}(\|H\|)
\end{array}\right]\quad {\rm and}\quad U_{\beta\beta}=Q+{O}(\|H\|^2), 
\end{equation}
Therefore, by noting that $\Pi_{\mathcal{S}_{+}^{\lvert\beta\rvert}}(\Lambda(Z)_{\beta \beta})={O}(\|H\|)$ and ${O}(\|H\|)$ is uniform for $A\in\mathbb{B}_r(\overline{A})$ with $\pi(\overline{A})=\pi(A)$, we obtain from the above discussion that 
\begin{align*}
R(A+H)-R(A)-R'(A;H)
&=\left[\begin{array}{ccc}
0 & 0 & 0 \\
0 & Q\Pi_{\mathcal{S}_{+}^{\lvert\beta\rvert}}(\Lambda(Z)_{\beta \beta})Q^T-\Pi_{\mathcal{S}_{+}^{\lvert\beta\rvert}}(H_{\beta \beta}) & 0 \\
0 & 0 & 0
\end{array}\right]+{O}(\|H\|^2)
\end{align*}
By \eqref{eq:unisocpi2} and \eqref{eq:unisocpi1}, we know that 
 \[
 \Lambda(Z)_{\beta\beta}=U_{\beta}^{T}\Lambda(A)U_{\beta}+U_{\beta}^{T}HU_{\beta}= U_{\beta\beta}^{T}H_{\beta\beta}U_{\beta\beta}+{O}(\|H\|^{2})= Q^{T}H_{\beta\beta}Q+{O}(\|H\|^{2})\,.
 \] Since $Q\in{\cal O}^{\lvert\beta\rvert}$, we have
 \[
 H_{\beta\beta}=Q\Lambda(Z)_{\beta\beta}Q^{T}+{O}(\|H\|^{2})\,,
 \]
 where ${O}(\|H\|)$ is uniform for $A\in\mathbb{B}_r(\overline{A})$ with $\pi(\overline{A})=\pi(A)$. 
Combining this with the globally Lipschitz continuity of $\Pi_{{\cal S}_{+}^{\lvert\beta\rvert}}(\cdot)$ and $\Pi_{{\cal S}_{+}^{\lvert\beta\rvert}}(Q \Lambda(Z)_{\beta\beta}Q^{T})=Q\Pi_{{\cal S}_{+}^{\lvert\beta\rvert}}( \Lambda(Z)_{\beta\beta})Q^{T}$, we obtain that
 \[
 Q\Pi_{{\cal S}_{+}^{\lvert\beta\rvert}}( \Lambda(Z)_{\beta\beta})Q^{T}-\Pi_{{\cal S}_{+}^{\lvert\beta\rvert}}(H_{\beta\beta})
 ={O}(\|H\|^{2})\,.
 \]
Therefore,
 \begin{equation}\label{eq:nonsmoothpart}
R(A+H)-R(A)-R'(A;H)={O}(\|H\|^{2})\,.
 \end{equation}
By combining (\ref{eq:unisocpi3}) and (\ref{eq:nonsmoothpart}), we know that for any ${\cal S}^{n}\ni H\to 0$,
\begin{equation}\label{eq:diag-strong-B-diff}
\Pi_{{\cal S}_{+}^{n}}(\Lambda(A)+H)-\Pi_{{\cal S}_{+}^{n}}(\Lambda(A))-\Pi'_{{\cal S}_{+}^{n}}(\Lambda(A);H)={O}(\|H\|^{2})\,
\end{equation}
 and ${O}(\|H\|^2)$ is uniform for $A\in\mathbb{B}_r(\overline{A})$ with $\pi(\overline{A})=\pi(A)$. 
Next, we consider the case that $A={P}^{T}\Lambda(A){P}$. Re-write (\ref{eq:unisocpi2}) as
\[
\Lambda(A)+{P}^{T}H{P}={P}^{T}U\Lambda(Z)U^{T}{P}\,.
\] Let $\widetilde{H}:={P}^{T}H{P}$. Then, we have 
$
\Pi_{{\cal S}^{n}_{+}}(A+H)={P}\,\Pi_{{\cal S}^{n}_{+}}(\Lambda(A)+\widetilde{H}){P}^{T}\,.
$
Therefore, since ${P}\in{\cal O}^{n}$, we know from (\ref{eq:diag-strong-B-diff}) and (\ref{eq:dd-projection}) that for any ${\cal S}^{n}\ni H\to 0$, (\ref{eq:piuni}) holds.

\section{Proof of Lemma \ref{lemma:unisec}} \label{app:C} 
We first consider the case where $A$ is diagonal. For notational simplicity, let $\Lambda=\Lambda(A)$ and $\Xi=\Lambda({A}+H)$. From \eqref{eq:lemap1}, we have $AU+HU=U\Xi$, 
which implies 
$$\Lambda_{\bar{\iota}_s\bar{\iota}_s}U_{\bar{\iota}_s\bar{\iota}_t}+(HU)_{\bar{\iota}_s\bar{\iota}_t}=U_{\bar{\iota}_s\bar{\iota}_t}\Xi_{\bar{\iota}_s\bar{\iota}_t}.$$
It follows that 
$$\Lambda_{\bar{\iota}_s\bar{\iota}_s}U_{\bar{\iota}_s\bar{\iota}_t}+\sum_{j=1}^{\bar{d}}H_{\bar{\iota}_s\bar{\iota}_j}U_{\bar{\iota}_j\bar{\iota}_t}=U_{\bar{\iota}_s\bar{\iota}_t}\Xi_{\bar{\iota}_t\bar{\iota}_t}.$$
This, together with Lemma \ref{lemma:unih} shows that 
\begin{align}
U_{\bar{\iota}_s\bar{\iota}_t}&=\Sigma^{st}\circ\sum_{j=1}^{\bar{d}}H_{\bar{\iota}_s\bar{\iota}_j}U_{\bar{\iota}_j\bar{\iota}_t}
=\Sigma^{st}\circ H_{\bar{\iota}_s\bar{\iota}_t}Q_t+{O}(\|H\|^2)\label{eq:lemapr2}
\end{align}
where $(\Sigma^{st})_{ij}=1/((\Xi_{\bar{\iota}_t\bar{\iota}_t})_i-(\Lambda_{\bar{\iota}_s\bar{\iota}_s})_j)$.  It is easy to see that $1/\big((\Xi_{\bar{\iota}_t\bar{\iota}_t})_i-(\Lambda_{\bar{\iota}_s\bar{\iota}_s})_j\big)=1/\big((\Lambda_{\bar{\iota}_t\bar{\iota}_t})_i-(\Lambda_{\bar{\iota}_s\bar{\iota}_s})_j\big)+{O}(\|H\|)$. Combining this with \eqref{eq:lemapr2}, we have 
$$U_{\bar{\iota}_s\bar{\iota}_t}=\Theta^{st}\circ H_{\bar{\iota}_s\bar{\iota}_t}Q_t+{O}(\|H\|^2),\;{\rm with}\;O(\|H\|^2)\;{\rm uniform\;for\;all}\;A\in\mathbb{B}_r(\overline{A}).$$
Next we consider $A=P^T\Lambda(A)P^T$. Re-write (\ref{eq:lemap1}) as
\[
\Lambda(A)+{P}^{T}H{P}={P}^{T}U\Lambda(A+H)U^{T}{P}\,.
\] Let $\widetilde{H}:={P}^{T}H{P}$. Since $P$ is an orthogonal matrix, the following proof is the same as the diagonal case. 
Thus we have completed the proof. 

}

\end{appendices}

\end{document}